\documentclass[preprint,12pt]{elsarticle}

\pdfoutput=1

\usepackage[english]{babel}
\usepackage{graphics}
\usepackage{graphicx}
\usepackage{epsfig}
\usepackage{amssymb}
\usepackage{amsmath}
\usepackage{amsthm}
\usepackage{amsfonts}
\usepackage{mathrsfs}

 \newtheorem{corollary}{Corollary}
 \newtheorem{algorithm}{Algorithm}
 \newtheorem{definition}{Definition}

\newtheorem{theorem}{Theorem}
\newtheorem{proposition}{Proposition}
\newtheorem{example}{Example}
\newtheorem{remark}{Remark}

\def\vec{\mathop{\mathrm{vec}}}
\def\diag{\mathop{\mathrm{diag}}}

\begin{document}

\begin{frontmatter}

\title{A contour integral approach to the computation of invariant pairs}

\author[xlim]{Moulay Barkatou}
\ead{moulay.barkatou@unilim.fr}

\author[xlim]{Paola Boito}
\ead{paola.boito@unilim.fr}

\author[xlim]{Esteban Segura Ugalde}
\ead{esteban.segura@etu.unilim.fr}

\address[xlim]{Universit\'e de Limoges - CNRS, XLIM UMR 7252, 123 avenue Albert Thomas, 87060 Limoges Cedex, France.}
  
\begin{abstract}
We study some aspects of the invariant pair problem for matrix polynomials, as introduced by Betcke and Kressner \cite{betcke2011perturbation} and by Beyn and Th\"ummler \cite{beyn2009continuation}. Invariant pairs extend the notion of eigenvalue-eigenvector pairs, providing a counterpart of invariant subspaces for the nonlinear case. We compute formulations for the condition numbers and the backward error for invariant pairs and solvents. We then adapt the Sakurai-Sugiura moment method \cite{asakura2010numerical} to the computation of invariant pairs, including some classes of problems that have multiple eigenvalues. Numerical refinement via a variant of Newton's method is also studied. Furthermore, we investigate the relation between the matrix solvent problem and the triangularization of matrix polynomials.

\end{abstract}
  
  
\end{frontmatter} 
 
 \section{Introduction}
 
 Invariant pairs, introduced and analyzed in \cite{betcke2011perturbation}, \cite{beyn2009continuation} and \cite{szyld2013several}, are a generalization of eigenpairs for matrix polynomials. Let $P(\lambda)=\sum_{j=0}^{\ell}A_j\lambda ^j$ be an $n\times n$ matrix polynomial, and choose a positive integer $k$. Then the matrices $X,S$ of sizes $n\times k$ and $k\times k$, respectively, form an invariant pair of size $k$ for $P(\lambda)$ if 
 $$
 P(X,S):=\sum_{j=0}^{\ell}A_jXS ^j=0.
 $$
 Invariant pairs offer a unified theoretical perspective on the problem of computing several eigenvalue-eigenvector pairs for a given matrix polynomial. From a numerical point of view, moreover, the computation of an invariant pair tends to be more stable than the computation of single eigenpairs, particularly in the case of multiple or tightly clustered eigenvalues. Observe that the notion of invariant pairs can also be applied to more general nonlinear problems, although here we will limit our presentation to matrix polynomials.
 
How to compute invariant pairs? Beyn and Th\"ummler (\cite{beyn2009continuation}) adopt a continuation method of predictor-corrector type. Betcke and Kressner (\cite{betcke2011perturbation}), on the other hand, establish a correspondence between invariant pairs of a given polynomial and of its linearizations. Invariant pairs for $P(\lambda)$ are extracted from invariant pairs of a linearized form and then refined via Newton's method. 

The approach we take in this paper to compute invariant pairs is based on contour integrals. Being able to specify the contour $\Gamma$ allows us to select invariant pairs that have eigenvalues in a prescribed part of the complex plane. Contour integrals play an important role in the definition and computation of moments, which form a Hankel matrix pencil yielding the eigenvalues of the given matrix polynomial that belong to the prescribed contour. The use of Hankel pencils of moment matrices is widespread in several applications such as control theory, signal processing or shape reconstruction, but nonlinear eigenvalue-eigenvector problems can also be tackled through this approach, as suggested for instance in \cite{asakura2010numerical} and \cite{beyn2012integral}. 
 E. Polizzi's FEAST algorithm \cite{polizzi2009density} is also an interesting example of contour-integral based eigensolver applied to large-scale electronic structure computations.
 
 Here we adapt such methods to the computation of invariant pairs. This work studies in particular the  scalar moment method and its relation with the multiplicity structure of the eigenvalues, but it also explores the behavior of the block version. We seek to compute invariant pairs
while avoiding the linearization of $P(\lambda)$, which explains the choice of moment methods. 
We are motivated, in  part, by the problem of better understanding the invariant pair problem in a symbolic or symbolic-numeric framework, that is, computing invariant pairs exactly, or with arbitrary accuracy via an effective combination of symbolic and numerical techniques: this is one of the reasons why we are interested in the issue of eigenvalue multiplicity.   

The last part of the paper shows an application of our results on invariant pairs to the particular case of matrix solvents, that is, to the matrix equation 
$$
 P(S):=\sum_{j=0}^{\ell}A_jS ^j=0.
$$

The matrix solvent problem has received remarkable attention in the literature, since Sylvester's work \cite{sylvester1884hamilton} in the 1880s. The relation between the Riccati and the  quadratic matrix equation is highlighted in \cite{bini2010transforming}, whereas a study on the  existence of solvents can be found in \citep{dennis1976algebraic}. Several works address the problem of computing a numerical approximation for the solution of the quadratic matrix equation:
an approach to compute, when possible, the dominant solvent is proposed in  \citep{dennis1978algorithms}. Newton's method and some variations are also used to approximate solvents numerically: see for example \cite{davis1981numerical}, \cite{kim2000numerical}, \cite{higham2001solving}, \cite{long2008improved}. The work in \cite{hashemi2010efficient} uses interval arithmetic to compute an interval matrix containing the exact solution to the quadratic matrix equation. 
For the case of the general matrix solvent problem, we can also cite \cite{bras1996spectral}, \cite{pereira2003solvents} and \cite{kratz1987numerical}.
On the other hand, the question of designing symbolic algorithms for computing solvents remains relatively unexplored. Attempts have been made to formulate the problem as a system of polynomial equations which can be solved via standard methods. However, this approach becomes cumbersome for problems of large size (see \cite{higham2000numerical}).

Here we exhibit computable formulations for the condition number and backward error of the general matrix solvent problem, generalizing existing works on the quadratic matrix equation. Moreover, we propose an adaptation of the moment method to the computation of solvents. Finally, we build on existing work on triangularization of matrix polynomials (see \cite{tisseur2013triangularizing} and \cite{taslaman2013triangularizing}) and explore the relationship between solvents of matrix polynomials in general and in triangularized form.


The paper is organized as follows. Section \ref{sec_preliminaries} introduces preliminary notions, definitions and notation. The backward error and condition number for the invariant pair problem are computed in Section \ref{sec_condition_ip}. 

Section \ref{sec_invpairs} is devoted to the computation of eigenvalues and invariant pairs through moments and Hankel pencils. Our main results here consist in Theorem \ref{theo_scalar}, Corollary \ref{coro_inv} and Theorem \ref{th_invariantpairs}. A comparison of different techniques for numerical refinement of invariant pairs is presented in Section \ref{sec_Newton}.

Finally, Sections \ref{sec_solvents}, \ref{sec_solvents_computation} and \ref{sec_solvents_triang} are devoted to the applications to matrix solvents.

The methods presented in the paper have been implemented both in a symbolic (exact) and in a numerical version. The Maple and Matlab codes
are available online at the URL\\ {\tt http://www.unilim.fr/pages\_perso/esteban.segura/software.html}.

\section{Matrix polynomials and invariant pairs}\label{sec_preliminaries}
A complex $n\times n$ {\em matrix polynomial} $P(\lambda)$ of degree $\ell$ takes the form: 
\begin{equation} 
P(\lambda) = A_0+A_1\lambda+A_2\lambda^2+ \cdots + A_\ell\lambda^\ell 
\end{equation}
where $A_\ell \neq 0$ and $A_i \in \mathbb{C}^{n\times n}$, for $i=0,\ldots,\ell$. \\
In this work, we assume that $P(\lambda)$ is regular, i.e., det$P(\lambda)$ does not vanish identically.

A crucial property of matrix polynomials is the existence of the Smith form (see, e.g., \cite{gohberg1982matrix}):
\begin{theorem}{[Thm. S1.1, \cite{gohberg1982matrix}]} Every $n\times n$  regular matrix polynomial $P(\lambda)$ admits the representation 
 \begin{equation} \label{Smith}
 D(\lambda) = E(\lambda)P(\lambda)F(\lambda),
 \end{equation}
 where $D(\lambda) = \diag{(d_1(\lambda),\ldots,d_n(\lambda))}$ is a diagonal polynomial matrix with monic scalar polynomials $d_i(\lambda)$ such that $d_i(\lambda)$ is divisible by $d_{i-1}(\lambda)$; $E(\lambda)$ and $F(\lambda)$ are matrix polynomials of size $n\times n$ with constant nonzero determinants. 
 \end{theorem}

The {\em polynomial eigenvalue problem} (PEP) consists in determining right eigenvalue-eigenvector pairs $(\lambda,x)\in\mathbb{C}\times\mathbb{C}^n$, with $x\neq 0$, such that
$$
P(\lambda)x=0,
$$
or left eigenvalue-eigenvector pairs $(\lambda,y)\in\mathbb{C}\times\mathbb{C}^n$, with $y\neq 0$, such that
$$
y^* P(\lambda)=0.
$$
A particular case of special interest is the {\em quadratic eigenvalue problem} (QEP), where $\ell=2$. Typical applications of the QEP include the vibrational analysis of various physical systems. A considerable amount of work has been done on the theoretical and computational study of the QEP: see for instance \cite{tisseur2001quadratic}.

As for the linear case, there is a notion of {\em algebraic} and {\em geometric multiplicity} of eigenvalues of matrix polynomials. If $\lambda_0$ is an eigenvalue of $P(\lambda)$, the algebraic multiplicity of $\lambda_0$ is its multiplicity as root of det$ P(\lambda)$, whereas the geometric multiplicity of $\lambda_0$ is the dimension of the null space of the matrix $P(\lambda_0)$. 

{\em Invariant pairs}, introduced and analyzed in \cite{betcke2011perturbation} and \cite{beyn2009continuation}, are a generalization of the notion of eigenpair for matrix polynomials.  
\begin{definition} \label{InvPairDef}
A pair $(X,S) \in \mathbb{C}^{n\times k}\times \mathbb{C}^{k\times k}$, $X\neq 0$, is called an invariant pair for $P(\lambda)$ if it satisfies the relation:
\begin{equation} \label{InvPair}
P(X,S):=A_\ell XS^\ell +\cdots+A_2XS^2 + A_1XS + A_0X = 0,
\end{equation}
where $A_i \in \mathbb{C}^{n\times n}$, $i=0,\ldots,\ell$, and $k$ is an integer between $1$ and $\ell n$. 
\end{definition}

The following definitions proposed in \cite{betcke2011perturbation} and \cite{gohberg1982matrix} will be helpful for our work, for instance, to allow for rank deficiencies in $X$.
\begin{definition} 
A pair $(X,S) \in \mathbb{C}^{n \times k}\times \mathbb{C}^{k\times k}$ is called minimal if there is $m \in \mathbb{N}$ such that:
\begin{equation*}
V_m(X,S) := \left[ \begin{matrix}
        XS^{m-1}\\ \vdots\\ XS\\ X    
\end{matrix} \right]
\end{equation*}
has full column rank. The smallest such $m$ is called minimality index of $(X,S)$.
\end{definition}

\begin{definition}
An invariant pair $(X,S)$ for a regular matrix polynomial $P(\lambda)$ of degree $\ell$ is called simple if $(X,S)$ is minimal and the algebraic multiplicities of the eigenvalues of $S$ are identical to the algebraic multiplicities of the corresponding eigenvalues of $P(\lambda)$. 
\end{definition}

Invariant pairs are closely related to the theory of standard pairs presented in \cite{gohberg1982matrix}, and in particular to {\em Jordan pairs}. If $(X,S)$ is a simple invariant pair and $S$ is in Jordan form, then $(X,S)$ is a Jordan pair.  

As an example, consider the following quadratic matrix polynomial, discussed in \cite{tisseur2001quadratic}:
\begin{equation} \label{Ex1}
P(\lambda) = \lambda^2 \left[ \begin{matrix}  1&0&0\\0&1&0\\ 0&0&0 \end{matrix} \right] + \lambda \left[ \begin{matrix}  -2&0&1\\0&0&0\\ 0&0&0 \end{matrix} \right] +\left[ \begin{matrix}  1&0&0\\0&-1&0\\ 0&0&1 \end{matrix} \right]
\end{equation}
It has eigenvalues $\lambda_1 = 1$ with algebraic multiplicity 3 and $\lambda_2 = -1$ with algebraic multiplicity 1. A corresponding Jordan pair $(X,J)$ is given by:
\begin{equation*}
X = \left[ \begin{matrix}  0&0&1&0\\1&1&0&1\\ 0&0&0&0 \end{matrix} \right], \; J = \text{diag}\left( -1,1,\left[ \begin{matrix} 1&1\\0&1\end{matrix} \right] \right)
\end{equation*} 
The notion of invariant pairs offers a theoretical perspective and a numerically more stable approach to the task of computing several eigenpairs of a matrix polynomial. Indeed, this problem is typically ill-conditioned in presence of multiple or nearly multiple eigenvalues, whereas the corresponding invariant pair formulation may have better stability properties.

In particular, simple invariant pairs play an important role when using a linearization approach as in \cite{betcke2011perturbation}, and ensure local quadratic convergence of Newton's method, as shown in \cite{kressner2009block}; see also \cite{szyld2013several}. 

\subsection{Condition number and backward error of the invariant pair problem}\label{sec_condition_ip}
In the following sections, we present explicit formulas for the backward error and the condition number of an invariant pair $(X,S)$. We follow the ideas presented in the articles \cite{tisseur2000backward} and \cite{higham2001solving}, which give expressions for backward errors and condition numbers for the polynomial eigenvalue problem and for a solvent of the quadratic matrix equation.

\subsubsection{Condition number} \label{SecCondNumInvPa}
Let $(X, S)$ be an invariant pair for $P(\lambda)$, and consider the perturbed polynomial
$$
(P+\Delta P)(\lambda) = (A_0+\Delta A_0) + \lambda (A_1+\Delta A_1) + \cdots +\lambda^\ell (A_\ell+\Delta A_\ell),
$$
where $\Delta P(\lambda) = \displaystyle \sum_{i=0}^\ell \lambda^i \Delta A_i$ and $\Delta A_0,\ldots, \Delta A_\ell \in \mathbb{C}^{n\times n}$. Let $(\Delta X,\Delta S)$ be a perturbation of $(X, S)$ such that
$$
(P + \Delta P)(X+\Delta X, S+\Delta S) = 0.
$$
A normwise condition number of the invariant pair $(X,S)$ can be defined as:
\begin{small}
\begin{align} \label{defCondNum}
\kappa (X,S) = \limsup_{\epsilon\rightarrow 0} \left\lbrace \frac{1}{\epsilon} \frac{\left\| \left[ \begin{array}{c} \Delta X \\ \Delta S \end{array}\right] \right\|_F}{\left\| \left[ \begin{array}{c} X \\  S \end{array}\right] \right\|_F} : ( P+\Delta P)(X+\Delta X, S+\Delta S)=0, \right. \nonumber \\ \| \Delta A_i\|_F \leq \epsilon \alpha_i, i=0, \ldots, \ell \bigg\}
\end{align}
\end{small}
The $\alpha_i$ are nonnegative weights that provide flexibility in how the perturbations are measured. A common choice is $\alpha_i=\|A_i\|_F$; however, if some coefficients are to be left unperturbed, $\Delta A_i$ can be forced to zero by setting $\alpha_i =0$.\\
\begin{theorem} \label{condNumInvPair}
The normwise condition number of the simple invariant pair $(X,S)$ is given by:
\begin{equation} \label{CondNumGen}
\kappa (X,S) = \frac{\left|\left|  \left[ \begin{matrix} B_{X} & B_{S} \end{matrix} \right]^+ B_{A}\right|\right|_2 }{ \left\| \left[ \begin{array}{c} X \\ S \end{array}\right] \right\|_F},
\end{equation}
where
\begin{align*}
B_{X} =& \displaystyle\sum_{j = 0}^{\ell} \left(  \left[ (S^j)^T \otimes A_j \right] \right), \quad B_{S} = \displaystyle\sum_{j = 1}^{\ell}  \displaystyle\sum_{i=0}^{j-1}\left( (S^{j-i-1})^T\otimes A_j X S^i \right),\\
B_{A} =& \left[ \begin{matrix} \alpha_{\ell} (X S^\ell)^T\otimes I_n & \cdots & \alpha_{0}X^T \otimes I_n  \end{matrix} \right].
\end{align*}
\begin{proof}
By expanding the first constraint in \eqref{defCondNum} and keeping only the first order terms, we get:
\begin{equation} \label{devep}
\displaystyle\sum_{j = 0}^{\ell} \Delta A_j X S^j + \displaystyle\sum_{j = 0}^{\ell} A_j \Delta X S^j +  \displaystyle\sum_{j = 1}^{\ell}A_jX\mathbb{D}S^j(\Delta S) = O(\epsilon^2),
\end{equation}
where $\mathbb{D}S^j$ denotes the Fr\'echet derivative of the map $S\mapsto S^j$: \index{Invariant pair!Fr\'echet derivative}
\begin{equation*}
\mathbb{D}S^j: \Delta S \mapsto \displaystyle\sum_{i=0}^{j-1}S^i \Delta S S^{j-i-1}.
\end{equation*}
Using on equation \eqref{devep} the $\vec{}$ operator and its property (see \cite{golub2013matrix}, pp. 28):
\begin{equation} \label{vecOpe}
\vec{(AXB)} = (B^T\otimes A)\vec{(X)},
\end{equation}
we obtain:
\begin{scriptsize}
\begin{align*}
\bullet & \vec{(\Delta P(X,S))} = \vec{\left( \displaystyle\sum_{j = 0}^{\ell} \Delta A_j X S^j\right) } =  \displaystyle\sum_{j = 0}^{\ell} \vec{\left(  \Delta A_j X S^j\right) } = \displaystyle\sum_{j = 0}^{\ell} \left(  \left[ (X S^j)^T \otimes I_n \right] \vec{(  \Delta A_j ) }\right) = \\
& = \left[ \begin{matrix} \alpha_{\ell} (X S^\ell)^T\otimes I_n & \cdots & \alpha_{0}X^T \otimes I_n  \end{matrix} \right] \left[ \begin{matrix} \vec{(  \Delta A_\ell ) }/\alpha_{\ell} \\ \vdots \\ \vec{(  \Delta A_0 ) }/\alpha_{0} \end{matrix} \right] =: B_{A} \vec{(  \Delta A)}, \\
\bullet & \vec{(P(\Delta X,S))} = \vec{\left( \displaystyle\sum_{j = 0}^{\ell} A_j \Delta X S^j\right) } =  \displaystyle\sum_{j = 0}^{\ell} \vec{\left( A_j \Delta X S^j \right) } = \displaystyle\sum_{j = 0}^{\ell} \left(  \left[ (S^j)^T \otimes A_j \right] \right) \vec{(  \Delta X ) } = \\
& =: B_{X} \vec{( \Delta X )}, \\
\bullet & \vec{\left(\displaystyle\sum_{j = 1}^{\ell} A_jX\mathbb{D}S^j(\Delta S) \right)} = \vec{\left(\displaystyle\sum_{j = 1}^{\ell} A_jX \displaystyle\sum_{i=0}^{j-1}S^i \Delta S S^{j-i-1} \right)} = \displaystyle\sum_{j = 1}^{\ell}  \displaystyle\sum_{i=0}^{j-1}\vec{\left( A_j X S^i \Delta S S^{j-i-1} \right)}=\\
&= \displaystyle\sum_{j = 1}^{\ell}  \displaystyle\sum_{i=0}^{j-1}\left( (S^{j-i-1})^T\otimes A_j X S^i \right) \vec{(\Delta S)} =: B_{S} \vec{( \Delta S )}.
\end{align*}
\end{scriptsize}
Then, we have:
$$
\left[ \begin{matrix} B_{X} & B_{S} \end{matrix} \right] y = -B_{A}x + O(\epsilon^2),
$$
where
\begin{align*}
y = \left[ \begin{matrix} \vec{( \Delta X )}\\ \vec{( \Delta S )} \end{matrix} \right], \quad and \quad x = \left[ \begin{matrix} \vec{(  \Delta A_\ell ) }/\alpha_{\ell} \\ \vdots \\ \vec{(  \Delta A_0 ) }/\alpha_{0} \end{matrix} \right]
\end{align*}
and therefore
$$
\| y \|_2 = \left\| \left[\begin{matrix}  \vec{(\Delta X)}\\ \vec{(\Delta S)} \end{matrix} \right]\right\|_2 =  \left\| \left[ \begin{array}{c} \Delta X \\ \Delta S \end{array}\right] \right\|_F.
$$
So we have that the definition \eqref{defCondNum} is equivalent to the following
\begin{footnotesize}
$$
\limsup_{\epsilon\rightarrow 0} \left\lbrace \frac{1}{\epsilon} \frac{\| y \|_2}{ \left\| \left[ \begin{array}{c} X \\ S \end{array}\right] \right\|_F} : \left[ \begin{matrix} B_{X} & B_{S} \end{matrix} \right] y = -B_{A}x + O(\epsilon^2), \| x\|_2 \leq \epsilon \right\rbrace = \frac{ \left\| \left[ \begin{matrix} B_{X} & B_{S} \end{matrix} \right]^+ B_{A}\right\|_2 }{ \left\| \left[ \begin{array}{c} X \\ S \end{array}\right] \right\|_F},
$$
\end{footnotesize}
where the matrix $\left[ \begin{matrix} B_{X} & B_{S} \end{matrix} \right]$ has full rank if the invariant pair $(X, S)$ is simple (see [Thm. 7, \cite{betcke2011perturbation}]).
\end{proof}
\end{theorem}

In order to better illustrate Theorem \ref{condNumInvPair}, let us consider the particular case $k=1$.
When $k = 1$, invariant pairs $(X,S)$ coincide with eigenpairs $(x,\lambda)$. In this case, the matrices $B_{X}$, $B_{S}$ and $B_{A}$ in \eqref{CondNumGen} are:
\begin{footnotesize}
\begin{align*}
B_{X} =& \displaystyle\sum_{j = 0}^{\ell} \left(  \left[ (\lambda^j)^T \otimes A_j \right] \right) = \displaystyle\sum_{j = 0}^{\ell} \left( \lambda^j A_j \right) = P(\lambda), \\
B_{S} =& \displaystyle\sum_{j = 1}^{\ell}  \displaystyle\sum_{i=0}^{j-1}\left( (\lambda^{j-i-1})^T\otimes A_j x \lambda^i \right) = \displaystyle\sum_{j = 1}^{\ell}  \displaystyle\sum_{i=0}^{j-1}\left( \lambda^{j-1} A_j x \right) = P'(\lambda) x,\\
B_{A} =& \left[ \begin{matrix} \alpha_\ell \lambda^\ell x^T\otimes I_n & \alpha_{\ell-1}\lambda^{\ell-1} x^T\otimes I_n & \cdots & \alpha_0 x^T \otimes I_n  \end{matrix} \right]
\end{align*}
\end{footnotesize}
Note that:
\begin{small}
\begin{align*}
B_{A}x =&\left[ \begin{matrix} \alpha_{\ell}\lambda^\ell x^T\otimes I_n & \cdots & \alpha_{0}x^T \otimes I_n  \end{matrix} \right] \left[ \begin{matrix} \vec{(\Delta A_\ell)}/\alpha_{\ell} \\ \vdots \\ \vec{(\Delta A_0)}/\alpha_{0} \end{matrix} \right] =\\
=& \vec{( \lambda^\ell  \Delta A_\ell x + \cdots + \Delta A_0 x )} = \vec{(\Delta P(\lambda) x)} = \Delta P(\lambda) x.
\end{align*}
\end{small}
Therefore, we obtain:
\begin{small}
\begin{align*}
&\left[ \begin{matrix} B_{X} & B_{S} \end{matrix} \right] y = -B_{A}x + O(\epsilon^2) \Leftrightarrow \left[ \begin{matrix} P(\lambda) & P'(\lambda) x \end{matrix} \right]  \left[ \begin{matrix} \Delta x \\ \Delta \lambda \end{matrix} \right] = -\Delta P(\lambda) x + O(\epsilon^2) \\
&\Leftrightarrow P(\lambda)\Delta x + P'(\lambda) x \Delta \lambda + \Delta P(\lambda) x = O(\epsilon^2).
\end{align*}
\end{small}
The last equation is consistent with the first part of the computation of the condition number for a nonzero simple eigenvalue $\lambda$ of $P(\lambda)$ presented in [Thm. 5, \cite{tisseur2000backward}]. The second part differs, because here we are estimating $\left\|\left[\begin{array}{c}\Delta x\\ \Delta\lambda\end{array}\right]\right\|_F$, whereas classical condition numbers for eigenvalue problems typically take into account angles between left and right eigenvectors. Of course, it would also be interesting to formalize a similar approach for invariant pairs, based on angles between suitable matrix manifolds (such as partially developed in \cite{betcke2011perturbation}).

\subsubsection{Backward error for $P(X,S)$} \label{backErrorInv}
Let $\alpha_i$, for $i=0,\ldots, \ell$, be nonnegative weights as in Section \ref{SecCondNumInvPa}. The backward error of a computed solution $(\tilde{X},\tilde{S})\in \mathbb{C}^{n\times k}\times \mathbb{C}^{k\times k}$ to \eqref{InvPair} can be defined as:
\begin{equation} \label{BEIP1}
\eta(\tilde{X}, \tilde{S})  = \min \{\epsilon: (P+\Delta P)(\tilde{X}, \tilde{S}) = 0, \| \Delta A_i\|_F \leq \epsilon \alpha_i, i=0,\ldots, \ell \}
\end{equation}
By expanding the first constraint in \eqref{BEIP1} we get:
\begin{equation} \label{BEIP2}
-P(\tilde{X}, \tilde{S}) = \Delta A_\ell \tilde{X} \tilde{S}^\ell +\cdots + \Delta A_0\tilde{X}.
\end{equation}
Then, we have
\begin{small}
\begin{equation*}
-P(\tilde{X}, \tilde{S}) = \left[ \begin{matrix}  \alpha_\ell^{-1}\Delta A_\ell & \ldots & \alpha_1^{-1}\Delta A_1 & \alpha_0^{-1}\Delta A_0 \end{matrix} \right]  \left[ \begin{matrix}  \alpha_\ell \tilde{X} \tilde{S}^\ell\\ \vdots \\ \alpha_1 \tilde{X}\tilde{S} \\ \alpha_0 \tilde{X} \end{matrix} \right] 
\end{equation*}
\end{small}
Taking the Frobenius norm, we obtain the lower bound for the backward error:
\begin{equation*}
\eta(\tilde{X}, \tilde{S}) \geq \frac{\| P(\tilde{X}, \tilde{S}) \|_F}{(\alpha_\ell^2 \| \tilde{X} \tilde{S}^\ell \|_F^2+ \cdots + \alpha_1^2 \| \tilde{X} \tilde{S} \|_F^2 + \alpha_0^2 \| \tilde{X} \|_F^2 )^{1/2}}.
\end{equation*}\\
Consider again equation \eqref{BEIP2}. Using \eqref{vecOpe}, we obtain:
\begin{footnotesize}
\begin{align*}
-\vec{(P(\tilde{X}, \tilde{S}))} & = ((\tilde{X} \tilde{S}^\ell)^T\otimes I_n )\vec{(\Delta A_\ell)} + \cdots + (\tilde{X}^T\otimes I_n)\vec{(\Delta A_0)} \\
& = \left[ \begin{matrix}  \alpha_\ell (\tilde{X} \tilde{S}^\ell)^T\otimes I_n & \ldots & \alpha_0 \tilde{X}^T\otimes I_n \end{matrix} \right]  \left[ \begin{matrix} \vec{(\Delta A_\ell)}/\alpha_\ell \\ \vdots \\ \vec{(\Delta A_0)}/\alpha_0 \end{matrix} \right],
\end{align*}
\end{footnotesize}
which can be written as:
\begin{equation} \label{BEIP3}
Hz = r, \, \, H \in \mathbb{C}^{nk\times (\ell+1)n^2}
\end{equation}
Here we assume that $H$ is full rank , to guarantee that \eqref{BEIP3} has a solution (backward error is finite). Then the backward error is the minimum 2-norm solution to:
\begin{equation} \label{BEIP4}
\eta(\tilde{X}, \tilde{S}) = \| H^+r\|_2,
\end{equation}
where  $H^+$ denotes the Moore-Penrose pseudoinverse of $H^+$. \\
Eq. \eqref{BEIP4} yields an upper bound for $\eta(\tilde{X}, \tilde{S})$:
\begin{equation*}
\eta(\tilde{X}, \tilde{S}) \leq \| H^+\|_2\| r\|_2 = \frac{\| r\|_2}{\sigma_{\min}(H)},
\end{equation*}
where $\sigma_{\min}$ denotes the smallest singular value, which is nonzero  by assumption. 
Note that:
\begin{small}
\begin{align*}
\sigma_{\min}(H)^2 & = \lambda_{\min} (HH^*) = \lambda_{\min}(\alpha_\ell^2(\tilde{X}\tilde{S}^\ell)^T\overline{\tilde{X}\tilde{S}^\ell}\otimes I_n + \cdots + \alpha_0^2 \tilde{X}^T\overline{X}\otimes I_n) \geq\\
& \geq \alpha_\ell^2\sigma_{\min}(\tilde{X}\tilde{S}^\ell)^2 + \cdots + \alpha_1^2\sigma_{\min}(\tilde{X}\tilde{S})^2 + \alpha_0^2\sigma_{\min}(\tilde{X})^2.
\end{align*}
\end{small}
Thus we obtain the upper bound for $\eta(\tilde{X}, \tilde{S})$:
\begin{equation*}
\eta(\tilde{X}, \tilde{S}) \leq \frac{\| P(\tilde{X}, \tilde{S})\|_F}{(\alpha_\ell^2\sigma_{\min}(\tilde{X}\tilde{S}^\ell)^2 + \cdots + \alpha_1^2\sigma_{\min}(\tilde{X}\tilde{S})^2 + \alpha_0^2\sigma_{\min}(\tilde{X})^2)^{1/2}}.
\end{equation*}

In the particular case $k = 1$, the approximate invariant pair $(\tilde{X}, \tilde{S})$ coincides with an approximate eigenpair $(\tilde{x},\tilde{\lambda})$. In this case, the definition \eqref{BEIP1} becomes: 
\begin{equation*}
\eta(\tilde{x},\tilde{\lambda})  = \min \{\epsilon: (P+\Delta P)(\tilde{x},\tilde{\lambda}) = 0, \| \Delta A_i\|_F \leq \epsilon \alpha_i, i=0,\ldots, \ell \},
\end{equation*}
which is the definition of the normwise backward error of an approximate eigenpair $(\tilde{x},\tilde{\lambda})$ for $P(\lambda)x=0$, presented in [(2.2), \cite{tisseur2000backward}].

\subsubsection{A numerical example}
Let us consider the {\em power plant problem} presented in \cite{betcke2013nlevp} and in \cite{tisseur2001quadratic}. This is a real symmetric QEP, with $P(\lambda)$ of size $8\times 8$, which describes the dynamic behaviour of a nuclear power plant simplified into an eight-degrees-of-freedom system. The problem is ill-conditioned due to the bad scaling of the matrix coefficients.\\
The maximum condition number for the eigenvalues of $P(\lambda)$, computed by the MATLAB function \texttt{polyeig}, is: 
$$ 
\kappa_{max} = \underset{\lambda \in \Lambda}{\max} \;\text{condeig}_\lambda =  1.0086\text{e+}008.
$$
Using the method that will be presented in Section \ref{sec_computingInvPairs} and Section \ref{sec_trapezoid}, we compute an invariant pair $(X,S)$ associated with the 11 eigenvalues with largest condition number inside the contour $\Gamma = \gamma + \rho e^{i\theta}$ ($\gamma =80+10i$, $\rho = 170$). The condition number and backward error for $(X,S)$ are
$$
\kappa (X,S) = 565.6746\qquad \textrm{and}\qquad
\eta(X,S) = 4.4548e-017.
$$
Observe that $\kappa (X,S)$ is significantly smaller than $ \kappa_{max}$.

\section{Computation of invariant pairs}\label{sec_invpairs}
Numerical methods based on contour integrals for the computation of eigenvalues of matrix polynomials and analytic matrix-valued functions have recently met with growing interest. Such techniques are related to the well-known method of moments, where the moments are computed by numerical quadrature.

In this section we explore a similar approach for computing invariant pairs. Our main reference is the Sakurai-Sugiura method (see \cite{asakura2010numerical} and \cite{sakurai2003projection}), as well as the presentation given in \cite{beyn2012integral}.

\subsection{The moment method and eigenvalues}
Let us begin by briefly recalling a few basic facts about the Sakurai-Sugiura moment method. Here we essentially follow the presentation given in \cite{asakura2010numerical}.

Let $\Gamma$ be a closed contour in the complex plane and let $u$ and $v$ be arbitrarily given vectors in $\mathbb{C}^n$. Define the function:
$$
f(\lambda):=u^H P(\lambda)^{-1} v.
$$
In the following, it will be understood that no eigenvalue of $P(\lambda)$ should lie exactly on the contour $\Gamma$: each eigenvalue should be either inside or outside the contour.

The next theorem, which can be found in \cite{asakura2010numerical}, gives a representation for $f(\lambda)$ that will prove useful later on. 

\begin{theorem}{[Thm. 3.1, \cite{asakura2010numerical}]}\label{form_of_f}
Let $D(\lambda) = \diag{(d_1(\lambda),\ldots,d_n(\lambda))}$ be the Smith form of $P(\lambda)$, and let $E(\lambda)$ and $F(\lambda)$ be as in \eqref{Smith}. Let $\chi_j(\lambda) = u^H \textbf{q}_j(\lambda) \textbf{p}_j(\lambda)^H v$, $1\leq j \leq n$. Then 
\begin{equation}\label{represent}
f(\lambda) = \displaystyle \sum_{j=1}^n \frac{\chi_j(\lambda)}{d_j(\lambda)},
\end{equation}
where $\textbf{q}_j(\lambda)$ and $\textbf{p}_j(\lambda)$ are the column vectors of $E(\lambda)$ and $F(\lambda)^H$, respectively. 
\end{theorem}

\begin{definition}
Let $k \in \mathbb{N}$. The $k$-th moment of $f(z)$ is:
\begin{equation} \label{moments}
\mu_k = \frac{1}{2\pi  \imath} \oint_\Gamma z^k f(z) dz.
\end{equation}
\end{definition}

For a positive integer $m$, define the Hankel matrices $H_0, H_1 \in \mathbb{C}^{m\times m}$ as follows:
 \begin{equation} \label{hankelMats}
 H_0 = \left[ \begin{matrix} \mu_0 & \mu_1 & \cdots & \mu_{m-1} \\ \mu_1 & \mu_2 & \cdots & \mu_{m} \\ \vdots & \vdots & & \vdots\\ \mu_{m-1} & \mu_m & \cdots & \mu_{2m-2} \end{matrix} \right], \quad  H_1 = \left[ \begin{matrix} \mu_1 & \mu_2 & \cdots & \mu_{m} \\ \mu_2 & \mu_3 & \cdots & \mu_{m+1} \\ \vdots & \vdots & & \vdots\\ \mu_{m} & \mu_{m+1} & \cdots & \mu_{2m-1} \end{matrix} \right]
 \end{equation}
 
 The eigenvalue algorithm presented in \cite{asakura2010numerical} relies on the following result: 
  \begin{theorem}{[Thm. 3.4, \cite{asakura2010numerical}]}\label{ss_eigenvalues}
Suppose that the polynomial $P(\lambda)$ has exactly $m$ eigenvalues $\lambda_1,\ldots,\lambda_m$ in the interior of $\Gamma$, and that these eigenvalues are distinct, simple and non degenerated. If $\chi_n(\lambda_\ell) \neq 0$ for $1\leq \ell \leq m$, then the eigenvalues of the pencil $H_1-\lambda H_0$ are given by $\lambda_1,\ldots,\lambda_m$.
 \end{theorem}
We now wish to investigate the behavior of the above method when the hypothesis that the $\lambda_i$'s are distinct is removed. In particular, we aim to generalize Theorem \ref{ss_eigenvalues}, which is based on the Vandermonde factorization of $H_0$ and $H_1$.

\begin{theorem} \label{theo_scalar}
Suppose that $P(\lambda)$ has exactly $m$ eigenvalues in the interior of $\Gamma$, namely, distinct eigenvalues $\lambda_0,\ldots,\lambda_s$ with algebraic multiplicities $m_0$, $\ldots$, $m_s$, respectively, such that $m=m_0+\cdots+m_s$. Moreover, assume that no eigenvalue of $P(\lambda)$ lies exactly on the contour $\Gamma$. If the geometric multiplicity of the $\lambda_i$'s, for $i=0,\ldots,s$, is equal to one, then the matrix $H_0$ is nonsingular and eigenvalues of the pencil $H_1-\lambda H_0$ are given by $\lambda_0,\ldots,\lambda_s$ with algebraic multiplicities $m_0, \ldots, m_s$.
\begin{proof} \quad 


Suppose that, in the Smith form \eqref{Smith} of $P(\lambda)$, the matrix $D(\lambda)$ is in the form $D(\lambda) = \text{diag}(d_1(\lambda),\ldots,d_{n-1}(\lambda),d_n(\lambda))$, where
$$
d_n(\lambda) = (\lambda-\lambda_0)^{m_0} (\lambda-\lambda_1)^{m_1}\cdots (\lambda-\lambda_s)^{m_s} \prod_{i=s+1}^{r}(\lambda-\lambda_i)^{m_i},
$$
and $\lambda_{s+1},\ldots,\lambda_r$ are the eigenvalues of $P(\lambda)$ located outside the contour $\Gamma$, with algebraic multiplicities $m_{s+1},\ldots,m_r$. Moreover, define $\tilde{d}_n(\lambda) = \prod_{i=0}^{s}(\lambda-\lambda_i)^{m_i}$, i.e., $\tilde{d}_n(\lambda)$ is the factor of $d_n(\lambda)$ whose roots are the eigenvalues of $P(\lambda)$ located inside $\Gamma$.

Since the geometric multiplicities of the $\lambda_i$'s are all equal to one, the factors $(\lambda - \lambda_i)$, for $i=0,\ldots,s$, do not appear in the monic scalar polynomials $d_0(\lambda), \ldots, d_{n-1}(\lambda)$.

Applying Theorem \ref{form_of_f}, we have:
\begin{align*}
\mu_k &= \frac{1}{2\pi  \imath} \oint_\Gamma z^k f(z) dz = \frac{1}{2\pi  \imath} \oint_\Gamma \displaystyle \sum_{j=1}^n \frac{\chi_j(z)}{d_j(z)}  z^k dz = \\
&= \frac{1}{2\pi  \imath} \oint_\Gamma \frac{\varphi (z)}{d_n(z)} z^k dz,
\end{align*}
where 
$$
\varphi (z) = \displaystyle \sum_{j=1}^{n} \chi_j(z) h_j(z), \; \text{with}\; h_j(z) = \frac{d_n(z)}{d_j(z)}.
$$
We can introduce partial fraction decompositions and write
\begin{align*}
\mu_k &= \frac{1}{2\pi  \imath} \oint_\Gamma \frac{\varphi (z)}{d_n(z)} z^k dz =\\
&=\frac{1}{2\pi  \imath} \oint_\Gamma \left( \displaystyle\sum_{i=1}^{m_0} \frac{c_{0,i} z^k}{(z-\lambda_0)^{i}}+\cdots+ \displaystyle\sum_{i=1}^{m_s} \frac{c_{s,i} z^k}{(z-\lambda_s)^{i}} \right) dz = \\
&= \displaystyle \sum_{j=0}^s \displaystyle \sum_{i=1}^{m_j} \frac{1}{2\pi  \imath} \oint_\Gamma \frac{c_{j,i} z^k}{(z-\lambda_j)^{i}} dz,
\end{align*}
where $c_{j,i }\in \mathbb{C}$, for $j=0,\ldots,s$ and $i=1,\ldots,m_j$. Classical results on residues then yield
\begin{align}
\mu_k &=  \displaystyle \sum_{j=0}^s \displaystyle \sum_{i=1}^{m_j}  c_{j,i} \text{Res}\left( \frac{z^k}{(z-\lambda_j)^i},\lambda_j \right) =\nonumber \\
&= \displaystyle \sum_{j=0}^s \displaystyle \sum_{i=1}^{m_j}  c_{j,i} \frac{1}{(i-1)!} \lim_{z \rightarrow \lambda_j} \frac{d^{i-1}}{dz^{i-1}}\left( (z-\lambda_j)^{i} \frac{z^k}{(z-\lambda_j)^{i}} \right) =\nonumber  \\
&= \displaystyle \sum_{j=0}^s \displaystyle \sum_{i=1}^{m_j}  c_{j,i} \frac{1}{(i-1)!} \lim_{z \rightarrow \lambda_j} \frac{d^{i-1}}{dz^{i-1}} \left( z^k \right) =\nonumber  \\
&= \displaystyle \sum_{j=0}^s \displaystyle \sum_{i=1}^{m_j} \nu_{j,i } \lambda_j^{k-i+1},\label{moment_formula}
\end{align}
where
$$
\nu_{j,i}=\left\{\begin{array}{l}
 \frac{c_{j,i}}{(i-1)!} (k-i+2)(k-i+3)\cdots k \quad\quad  {\rm if }\quad\quad k\geq i-1,\\
 0\quad\quad  {\rm otherwise. }
 \end{array}\right.
$$
Now, consider the pencil $H_1-\lambda H_0$ with $H_0$ and $H_1$ defined as in \eqref{hankelMats}.
Because of \eqref{moment_formula}, and of the fact that $\lambda_0,\ldots,\lambda_s$ are roots of $\tilde{d}_n(\lambda)$, the moments $\mu_k$ satisfy  a linear recurrence equation of the form:
\begin{equation} \label{recurrence}
\mu_k = a_{m-1}\mu_{k-1} + a_{m-2}\mu_{k-2}+\cdots + a_0 \mu_{k-m}.
\end{equation}
Moreover, $\tilde{d}_n(\lambda)$ is the polynomial of smallest degree that has roots $\lambda_0,\ldots,\lambda_s$ with the prescribed multiplicities $m_0,\ldots, m_s$, so the recurrence \eqref{recurrence} has the shortest possible length; also, note that the $a_i$'s in \eqref{recurrence} are actually the coefficients of $\tilde{d}_n(\lambda)$.
Therefore, the matrices $H_0$ and $H_1$ have full rank. The same argument shows that $H_0$ and $H_1$ are rank-deficient if taken of size larger than $m\times m$ (see \cite{meurant2010matrices}, [\cite{gantmacher1960theory}, Vol. 2, pp. 205)] and \cite{boley1997general}).


As a consequence of the shifted Hankel form of $H_0$ and $H_1$,  we have 
$H_0 C = H_1$, where $C$ is a matrix in companion form
\begin{equation*}
C = \left[ \begin{matrix} 0 & 0 & \cdots & 0& x_0 \\ 1 & 0 & \cdots& 0 & x_1 \\ 0&1&\cdots&0&x_2\\ \vdots & \vdots &\ddots&\vdots & \vdots\\  0 & 0 & \cdots & 1&x_{m-1} \end{matrix} \right],
\end{equation*}
and its last column is given by the solution of the linear system
\begin{equation} \label{system}
H_0 \left[ \begin{matrix} x_0\\ x_1\\ \vdots\\ x_{m-1}\end{matrix} \right] = \left[ \begin{matrix} \mu_m\\ \mu_{m+1}\\ \vdots\\ \mu_{2m-1}\end{matrix} \right].
\end{equation}
The polynomial of degree $m$:
\begin{equation*}
p(\lambda) = \lambda^m -x_{m-1}\lambda^{m-1} - \cdots - x_0
\end{equation*}
is a scalar multiple of $(\lambda-\lambda_0)^{m_0} (\lambda-\lambda_1)^{m_1}\cdots (\lambda-\lambda_s)^{m_s}$, and its roots are the $\lambda_i$'s generating the entries of the pencil $H_1-\lambda H_0$. So we also have that the $\mu_i$'s satisfy the recurrence \eqref{recurrence}:
\begin{equation*}
\mu_k = x_{m-1}\mu_{k-1} + x_{m-2}\mu_{k-2}+\cdots + x_0 \mu_{k-m},
\end{equation*}
where $k=m,m+1,\ldots,2m$.

Consider now the Jordan matrix
$$
J = \begin{bmatrix} J_1 & \;     & \; \\ \;  & \ddots & \; \\ \;  & \;     & J_s\end{bmatrix},
$$
where each block $J_i$, of dimension $m_i$, is a square matrix of the form
$$
J_i =  \begin{bmatrix} \lambda_i & 1  & \;   & \;  \\ \;   & \lambda_i    & \ddots & \;  \\ \;   & \;   & \ddots & 1  \\ \;   & \;   & \;     & \lambda_i  \end{bmatrix},
$$
and define the confluent Vandermonde matrix 
$$
V = \left( \begin{matrix} \textbf{v} & J^T\textbf{v} & \cdots & (J^T)^{r-1}\textbf{v} \end{matrix} \right),
$$
where $\textbf{v}^T = \left( \begin{matrix} e_1^{[m_1]T} & \cdots & e_1^{[m_s]T} \end{matrix} \right)$ is  partitioned conformally with $J$ and $e_1^{[m_\ell]T} = \left( \begin{matrix} 1 & 0 & \cdots & 0 \end{matrix} \right)^T$ is the $m_\ell$-dimensional unit coordinate vector. Then we have
\begin{align*}
V C &= \left( \begin{matrix} \textbf{v} & J^T\textbf{v} & \cdots & (J^T)^{r-1}\textbf{v} \end{matrix} \right) C\\
&= \left( \begin{matrix} J^T\textbf{v} & \cdots & (J^T)^{r-1}\textbf{v} & -(x_0 I+x_1J+\cdots+x_{m-1}J^{r-1})^T\textbf{v} \end{matrix} \right)\\
&= \left( \begin{matrix} J^T\textbf{v} & \cdots & (J^T)^{r-1}\textbf{v}&(J^T)^{r}\textbf{v} \end{matrix} \right) \\
&= J^TV,
\end{align*} 
where we have used the property $p(J)=0$, that is, the Cayley-Hamilton theorem. 

We can now introduce the Vandermonde decomposition of the Hankel matrices $H_0$ and $H_1$. 
From the results presented in \cite{boley1997general} it follows that 
there exist block matrices $B_0 = \text{diag}(D_1^{(0)},\ldots,D_s^{(0)})$ and $B_1 = \text{diag}(D_1^{(1)},\ldots,D_s^{(1)})$, partitioned conformally with $J$, satisfying the conditions $B_0J^T = JB_0$ and $B_1J^T = JB_1$,
so that
$$
H_i = V^TB_iV, \quad \text{for} \; i=0,1.
$$
Moreover, we can prove that $J B_0 = B_1$:
\begin{align*}
H_0 C = H_1&\Leftrightarrow (V^TB_0V) C = V^TB_1V \Leftrightarrow V^TB_0J^T V = V^TB_1V \\
&\Leftrightarrow  V^T J B_0 V = V^TB_1V \Leftrightarrow J B_0 = B_1,
\end{align*}
where we used the properties $ VC=J^TV$ and $ B_iJ^T = JB_i$.\\
Therefore, we have:
\begin{align*}
H_1 - \lambda H_0 &= V^TB_1V-\lambda V^TB_0V = V^TJ B_0 V-\lambda V^TB_0V \\
&=V^T(J-\lambda I)B_0V.
\end{align*}
So the eigenvalues of $H_1 - \lambda H_0$ are $\lambda_0,\ldots,\lambda_s$ with respective multiplicities $m_0,\ldots,m_s$.
\end{proof}
\end{theorem}

\begin{example}
 Consider the matrix polynomial:
\begin{equation*}
P(\lambda) = \lambda^2 I+ \lambda \left[ \begin{matrix} -2&0&0&0\\ 0&-1&0&0\\ 0&0&-6&0\\ 0&0&0&-5 \end{matrix} \right] + \left[ \begin{matrix} 1&0&0&1\\ 0&\frac{1}{4}&0&0\\ 0&0&9&0\\ 0&0&0&6 \end{matrix} \right].
\end{equation*}
$P(\lambda)$ has eigenvalues: $\lambda_1=\frac{1}{2}, \lambda_2 = 1, \lambda_3 = 2, \lambda_4 = 3$ with algebraic multiplicities: $m_1=2, m_2 = 2, m_3 = 1, m_4 = 3$. The associated Smith form is:
$D(\lambda) = \diag{(\left(1,1,1,(\lambda-\frac{1}{2})^2(\lambda-1)^2(\lambda-2)(\lambda-3)^3\right)}$.


Choosing vectors $u = \left[ \begin{matrix} 2&-2&1&-1 \end{matrix} \right]^T$ and  $v = \left[ \begin{matrix} 0&1&0&2 \end{matrix} \right]^T$, we find that:
\begin{equation*}
 H_0 = \left[ \begin{matrix} -3&-7&-9&\frac{-21}{2}\\ -7&-9&\frac{-21}{2}&-12\\ -9&\frac{-21}{2}&-12&\frac{-109}{8}\\ \frac{-21}{2}&-12&\frac{-109}{8}&\frac{-123}{8} \end{matrix} \right], \qquad  H_1 = \left[ \begin{matrix} -7&-9&\frac{-21}{2}&-12\\ -9&\frac{-21}{2}&-12&\frac{-109}{8}\\ \frac{-21}{2}&-12&\frac{-109}{8}&\frac{-123}{8}\\ -12&\frac{-109}{8}&\frac{-123}{8}&\frac{-551}{32} \end{matrix} \right].
\end{equation*}
Then, we have:
\begin{equation*}
C = H_0^{-1}H_1 =  \left[ \begin{matrix} 0&0&0&\frac{-1}{4}\\ 1&0&0&\frac{3}{2}\\ 0&1&0&\frac{-13}{4}\\ 0&0&1&3 \end{matrix} \right].
\end{equation*}
Note that the eigenvalues of $C$ are $\frac{1}{2},\frac{1}{2},1,1$. Moreover, the companion matrix $C$ is associated with the monic polynomial:
\begin{equation*} 
p(\lambda) = \lambda^4 - 3\lambda^3 + \frac{13}{4}\lambda^2 - \frac{3}{2}\lambda +\frac{1}{4},
\end{equation*}
whose roots are indeed $\frac{1}{2},\frac{1}{2},1,1$.

In fact, the Vandermonde factorization for $H_0$ and $H_1$ is $H_i=V^TB_iV$, $i=0,1$, where
\begin{footnotesize}
$$
V = \left[ \begin{matrix} 1&\frac{1}{2}&\frac{1}{4}&\frac{1}{8}\\ 0&1&1&\frac{3}{4}\\ 1&1&1&1\\ 0&1&2&3 \end{matrix} \right],\quad
B_0 = \left[ \begin{matrix} 0&-2&0&0\\-2&0&0&0\\ 0&0&-3&-2\\ 0&0&-2&0 \end{matrix} \right], \quad
B_1 = \left[ \begin{matrix} -2&-1&0&0\\-1&0&0&0\\ 0&0&-5&-2\\ 0&0&-2&0 \end{matrix} \right],
$$
\end{footnotesize}
and $JB_0=B_1$, with
$$
J=\left[ \begin{matrix} \frac{1}{2}&1&0&0\\ 0&\frac{1}{2}&0&0\\ 0&0&1&1\\ 0&0&0&1 \end{matrix} \right].
$$
\end{example}

\begin{example} \label{examp}
Consider the quadratic matrix polynomial:
$$
P(\lambda) = \lambda^2 \left[ \begin{matrix} 1&0&0\\ 2&1&0\\ -1&1&-2 \end{matrix} \right]  + \lambda \left[ \begin{matrix} 0&0&0\\ -4&-2&0\\ 2&-2&4 \end{matrix} \right] + \left[ \begin{matrix} -2&1& -2\\ 2&1&0\\ -1&1&-2 \end{matrix} \right].
$$
with associated Smith form:
$$
D(\lambda) =  \diag{(\left(d_1(\lambda), d_2(\lambda), d_3(\lambda)\right)} = \left[ \begin{matrix} 1&0&0\\ 0&(\lambda-1)^2&0\\ 0&0&(\lambda-1)^3(\lambda+1) \end{matrix} \right].
$$
The Jordan matrix associated with the linearized form of $P(\lambda)$ is:
$$
J = \begin{bmatrix} J_1 & 0\\ 0 & J_2\end{bmatrix},
$$
where each block $J_i$ is:
$$
J_1 =  \begin{bmatrix} 1 & 1\\ 0&1 \end{bmatrix}, \quad J_2 =  \begin{bmatrix} -1&0&0&0\\ 0&1&1&0\\ 0&0&1&1\\ 0&0&0&1 \end{bmatrix}.
$$
Note that we have the eigenvalue $\lambda = 1$ in different Jordan blocks.\\
Choose $\Gamma$ as the circle $\varphi(t) = 1+\frac{1}{10} e^{\imath t}$, which contains 5 eigenvalues $\lambda = 1$, and consider vectors $u = \left[ \begin{matrix} 3&1&-2 \end{matrix} \right]^T$ and  $v = \left[ \begin{matrix} 3&-1&-2 \end{matrix} \right]^T$. Theorem \ref{theo_scalar} implies that the moment method yields the eigenvalues of $P(\lambda)$ inside the contour, which are roots of $d_3(\lambda)$, i.e., $\lambda_0 = 1$ with multiplicity $m_0 = 3$. Let us compute the matrices $H_0$ and $H_1$ of size $3\times 3$:
$$
 H_0 = \left[ \begin{matrix} 7&3&3\\ 3&3&7\\ 3&7&15 \end{matrix} \right], \qquad  H_1 = \left[ \begin{matrix} 3&3&7\\ 3&7&15\\ 7&15&27 \end{matrix} \right].
$$
The matrix $H_0$ is nonsingular. Then, we have:
\begin{equation*}
C = H_0^{-1}H_1 =  \left[ \begin{matrix} 0&0&1\\ 1&0&-3\\ 0&1&3 \end{matrix} \right].
\end{equation*}
Note that eig$(C) = 1,1,1$.

As the contour contains 5 eigenvalues, we might ask what happens when taking the Hankel matrix $\hat{H}_0$ of size $5\times 5$:
$$
 \hat{H}_0 = \left[ \begin{matrix} 7&3&3&7&15\\ 3&3&7&15&27\\ 3&7&15&27&43\\ 7&15&27&43&63\\ 15&27&43&63&63 \end{matrix} \right].
$$
This matrix is singular, therefore we will not be able to find all the 5 eigenvalues inside the contour. The method miss the additional multiplicities associated with the polynomial $d_2(\lambda)$.
\end{example}

\begin{remark}
We conclude that the scalar moment method can be used to compute the (possibly multiple) eigenvalues of $P(\lambda)$ that belong to the interior of $\Gamma$ and that are roots of the polynomial $d_n(\lambda)$. The method misses the additional multiplicities associated with the polynomials $d_1(\lambda),\ldots,d_{n-1}(\lambda)$.
\end{remark}
The above remark is consistent with the fact that the Jordan form of a companion matrix only contains one Jordan block for each eigenvalue: it is not possible to capture multiple eigenvalues associated with several Jordan blocks.

In order to ``see'' the additional eigenvalues that are roots of $d_1(\lambda),\ldots,d_{n-1}(\lambda)$, the block version of the moment method may be useful (see Section \ref{block_mome}).

\subsection{Computing invariant pairs via moment pencils} \label{sec_computingInvPairs}
Let $\Gamma$ be a closed contour, $\lambda_1,\ldots,\lambda_m$ all the eigenvalues of $P(\lambda)$ in the interior of $\Gamma$ and the matrices $H_0$ and $H_1$ defined as in \eqref{hankelMats}.\\
For $k=0,1,\ldots,m-1$ and a nonzero vector $v \in \mathbb{C}^n$, consider the vectors:
\begin{equation}\label{eq_sk}
s_k = \frac{1}{2\pi  \imath} \oint_\Gamma z^k P(z)^{-1} v dz,
\end{equation}
The method proposed in \cite{asakura2010numerical} for the computation of the eigenvectors of $P(\lambda)$ is based on the following result.

\begin{theorem}{[Thm. 3.5, \cite{asakura2010numerical}]}\label{ss_eigenvectors}
Let $(\lambda_i,w_i)$, $i=1,\ldots,m$ be eigenpairs for the matrix pencil $H_1-\lambda H_0$, where the simple, distinct, nondegenerate eigenvalues $\lambda_i$ belong to the interior of a given closed contour $\Gamma$. Let $S=[s_0,\cdots ,s_{m-1}]$. Then, for $ i=1,\ldots,m$, the vector $y_i=Sw_i$ 
is an eigenvector of $P(\lambda)$ corresponding to the eigenvalue $\lambda_i$.
\end{theorem}
Theorem \ref{ss_eigenvectors} is readily applied to invariant pairs.
\begin{corollary} \label{coro_inv}
With the hypotheses of Theorem \ref{ss_eigenvectors}, $S = \left[ \begin{matrix} s_0,&s_1,&\ldots,&s_{m-1}  \end{matrix} \right]$ and $C=H_0^{-1}H_1$. Then the pair $(S,C)$ satisfies $P(S, C)=0$, i.e., $(S, C)$ is a simple invariant pair for $P(\lambda)$.
\begin{proof}
Note that the pair $(Y,\Lambda)$, where $\Lambda=\textnormal{diag}(\lambda_1,\ldots,\lambda_m)$ and $Y=[y_1,\ldots,y_m]$, is clearly an invariant pair for $P(\lambda)$, that is,
$$
P(Y,\Lambda)=\sum_{j=0}^{\ell}A_jY\Lambda^j=0.
$$
 Moreover, we know that $C=H_0^{-1}H_1=V^{-1}\Lambda V$, where $V$ is the classical Vandermonde matrix associated with $\lambda_1,\ldots,\lambda_m$, and that the columns of $V^{-1}$ are eigenvectors of $H_1-\lambda H_0$. So we have
\begin{align*}
0&=\sum_{j=0}^{\ell}A_jY\Lambda^j=\sum_{j=0}^{\ell}A_jY\Lambda^jV=\\
&=\sum_{j=0}^{\ell}A_jYVV^{-1}\Lambda^jV=\sum_{j=0}^{\ell}A_jSC^j=
P(S,C),
\end{align*}
that is, $(S,C)$ is also an invariant pair of $P(\lambda)$.
\end{proof}
\end{corollary}
What can we say about more general cases, where some of the hypotheses of Theorem \ref{ss_eigenvectors} are removed? If we remove the hypothesis that the $\lambda_i$'s are distinct, we can prove the following. 
\begin{theorem}\label{th_invariantpairs}
With the hypotheses of Theorem \ref{theo_scalar}, let $S = \left[ \begin{matrix} s_0,&s_1,&\ldots,&s_{m-1}  \end{matrix} \right]$ and $C=H_0^{-1}H_1$. Then the pair $(S,C)$ satisfies $P(S,C)=0$, i.e., $(S,C)$ is a simple invariant pair for $P(\lambda)$.
\begin{proof}
Consider again the columns $q_1(\lambda),\ldots, q_n(\lambda)$ of the matrix $F(\lambda)$ in the Smith form \eqref{Smith} and the definition of $s_k$ given in \eqref{eq_sk}. A similar computation to \eqref{moment_formula} shows that
$$
S=[s_0,\ldots,s_{m-1}]=\mathcal{Q}V,
$$
where
\begin{align*}
&\mathcal{Q}=[\mathcal{Q}_0,\ldots,\mathcal{Q}_s],\\
&\mathcal{Q}_j=[\gamma_{0,j}q_n(\lambda_j),\,\gamma_{1,j}q_n'(\lambda_j),\ldots, \gamma_{m_j-1,j}q_n^{(m_j-1)}(\lambda_j)], \,\, {\rm for }\, j=0,\ldots,s,
\end{align*}
the $\gamma_{i,j}$'s are complex coefficients and $V$ is the confluent Vandermonde matrix defined above.

It is shown in \cite{asakura2010numerical} (Lemma 2.4) that, if a complex number $\zeta$ is a root of $d_j(\lambda)$ for some index $1\leq j\leq n$, then $P(\zeta)q_j(\zeta)=0$.
In our case, this implies that the vector $q_n(\lambda)$ is a root polynomial of $P(\lambda)$ corresponding to the eigenvalue $\lambda_j$, for each $j=0,\ldots,s$; see \cite{gohberg1982matrix}, section 1.5, for the definition and properties of root polynomials. It follows that $[q_n(\lambda_j),\,q_n'(\lambda_j),\ldots,q_n^{(m_j-1)(\lambda_j)}]$ forms a Jordan chain for the eigenvalue $\lambda_j$. So we have that $(\mathcal{Q},J)$ is an invariant pair for $P(\lambda)$. Moreover, if $C=H_0^{-1}H_1$ as usual, we have
\begin{align*}
0&=\sum_{j=0}^{\ell}A_j\mathcal{Q}J^j=\sum_{j=0}^{\ell}A_j\mathcal{Q}J^jV=\\
&=\sum_{j=0}^{\ell}A_j\mathcal{Q}VV^{-1}J^jV=\sum_{j=0}^{\ell}A_jSC^j=
P(S,C),
\end{align*}
Therefore, $(S,C)$ is a simple invariant pair for $P(\lambda)$.
\end{proof}
\end{theorem}

\begin{example}
 Consider the matrix polynomial:
\begin{equation*}
P(\lambda) = \lambda^2 \left[ \begin{matrix}1&0\\ 0&1 \end{matrix} \right] + \lambda \left[ \begin{matrix} -2&0\\ 2&-1 \end{matrix} \right] + \left[ \begin{matrix} 1&0\\ 0&0 \end{matrix} \right],
\end{equation*}
which has eigenvalues $\lambda_1=0$ with algebraic multiplicity 1 and $\lambda_2 =1$ with algebraic multiplicity 3.\\
Suppose we are interested in the eigenvalues $\lambda_2$. Then we can choose a contour $\Gamma(t) = z_0 + R e^{\imath t}$, $t\in [0,2\pi]$, where $z_0=1$ and $R = \frac{1}{2}$. 

Choosing the vectors $u = \left[ \begin{matrix} 1 & -1 \end{matrix} \right]^T$ and  $v = \left[ \begin{matrix} -1 & 1 \end{matrix} \right]^T$, we find that:
\begin{equation*}
 H_0 = \left[ \begin{matrix} -1&-2&-5\\ -2&-5&-10\\ -5&-10&-17 \end{matrix} \right], \qquad  H_1 = \left[ \begin{matrix} -2&-5&-10\\ -5&-10&-17\\ -10&-17&-26 \end{matrix} \right]
\end{equation*}
Then, we have that the pair $(S,C)$ given by Theorem \ref{th_invariantpairs}
\begin{equation*}
S = \left[ \begin{matrix} 0&-1&-2\\ 1&1&3 \end{matrix} \right] \quad \text{and} \quad C = H_0^{-1}H_1 =  \left[ \begin{matrix} 0&0&1\\ 1&0&-3\\ 0&1&3 \end{matrix} \right]
\end{equation*}
is an invariant pair, i.e., it satisfies $P(S,C)=0$.\\
Note that the companion matrix $C$ is associated with the monic polynomial:
\begin{equation*} 
p(\lambda) = \lambda^3 - 3\lambda^2 + 3\lambda -1,
\end{equation*}
which has as roots: 1,1,1.
\end{example}
 
 \subsection{The block moment method} \label{block_mome}
 Instead of the scalar version of the moment method, we can consider a Hankel pencil constructed by block moments $M_k \in \mathbb{C}^{\xi\times \xi}$, for a suitable positive integer $\xi$.
\begin{definition}
Let $k$ be a positive integer and $U, V \in \mathbb{C}^{n \times \xi}$ nonzero matrices with linearly independent columns. For $k=0,1,\ldots$, define the block moment $M_k \in \mathbb{C}^{\xi\times \xi}$ as:
$$
M_k = \frac{1}{2\pi  \imath} \oint_\Gamma z^k U^H P(z)^{-1} V dz.
$$
\end{definition}
Then the block Hankel matrices $H_{\xi0}, H_{\xi1} \in \mathbb{C}^{\tilde{m}\xi\times \tilde{m}\xi}$ are defined as:
 \begin{equation*}
 H_{\xi0} = \left[ \begin{matrix} M_0 & M_1 & \cdots & M_{\tilde{m}-1} \\ M_1 & M_2 & \cdots & M_{\tilde{m}} \\ \vdots & \vdots & & \vdots\\ M_{\tilde{m}-1} & M_{\tilde{m}} & \cdots & M_{2\tilde{m}-2} \end{matrix} \right], \quad  H_{\xi1} = \left[ \begin{matrix} M_1 & M_2 & \cdots & M_{\tilde{m}} \\ M_2 & M_3 & \cdots & M_{\tilde{m}+1} \\ \vdots & \vdots & & \vdots\\ M_{\tilde{m}} & M_{\tilde{m}+1} & \cdots & M_{2\tilde{m}-1} \end{matrix} \right]
 \end{equation*}
 Polynomial eigenvalue computation via the eigenvalues of the pencil $H_{\xi1}-\lambda H_{\xi0}$ is discussed in \cite{asakura2010numerical} and \cite{beyn2012integral}. See also \cite{leblanc2013solving} for an application to acoustic nonlinear eigenvalue problems.\\
Invariant pairs can be computed from block moments by applying an approach that is similar to the one described in the previous section for the scalar version. For $k=0,1,\ldots, \tilde{m}-1$, consider the matrices $S_k \in \mathbb{C}^{n\times \xi}$ defined as:
$$
S_k=\frac{1}{2\pi  \imath} \oint_\Gamma z^k  P(z)^{-1} V dz.
$$
Then, we have the following result.\\
\begin{proposition}
Let $\Gamma$ be a closed contour, let the block Hankel matrix $H_{\xi0} \in \mathbb{C}^{\tilde{m}\xi\times \tilde{m}\xi}$ be nonsingular and $m$ be the number of eigenvalues inside of $\Gamma$. If $\tilde{m}\xi=m$ and $Y=[S_0,\ldots,S_{\tilde{m}-1}]$, $T=H_{\xi0}^{-1}H_{\xi1}$, then the pair $(Y,T)$ satisfies $P(Y,T)=0$, i.e., $(Y,T)$ is a simple invariant pair for $P(\lambda)$.
\end{proposition}
With the condition that the size of the block Hankel matrix $H_{\xi0} \in \mathbb{C}^{\tilde{m}\xi\times \tilde{m}\xi}$ is equal to the number of eigenvalues inside of $\Gamma$, i.e., if $\tilde{m}\xi=m$, we get:
\begin{equation*}
T = H_{\xi0}^{-1}H_{\xi1} =  \left[ \begin{matrix} 0 & 0 & \cdots & 0& -X_0 \\ I & 0 & \cdots& 0 & -X_1 \\ 0&I&\cdots&0&-X_2\\ \vdots & \vdots &\ddots&\vdots & \vdots\\  0 & 0 & \cdots & I&-X_{m-1} \end{matrix} \right],
\end{equation*}
where
\begin{equation*}
\left[ \begin{matrix} -X_0\\ -X_1\\ \vdots\\ -X_{m-1}\end{matrix} \right] = H_{\xi0}^{-1} \left[ \begin{matrix} M_m\\ M_{m+1}\\ \vdots\\ M_{2m-1}\end{matrix} \right].
\end{equation*}
Consequently, since $T$ has a block companion form, the problem of finding the eigenvalues $\lambda_1,\ldots,\lambda_m$ of  is equivalent to the problem of finding the eigenvalues of the matrix polynomial:
$$
L(\lambda) := \lambda^{\ell}+X_{\ell-1}\lambda^{\ell-1}+\cdots+X_{1}\lambda+X_{0} = 0.
$$
\begin{example}
Consider again the matrix polynomial of Example \ref{examp}:
$$
P(\lambda) = \lambda^2 \left[ \begin{matrix} 1&0&0\\ 2&1&0\\ -1&1&-2 \end{matrix} \right]  + \lambda \left[ \begin{matrix} 0&0&0\\ -4&-2&0\\ 2&-2&4 \end{matrix} \right] + \left[ \begin{matrix} -2&1& -2\\ 2&1&0\\ -1&1&-2 \end{matrix} \right].
$$
with associated Smith form:
$$
D(\lambda) =  \diag{(\left(d_1(\lambda), d_2(\lambda), d_3(\lambda)\right)} = \left[ \begin{matrix} 1&0&0\\ 0&(\lambda-1)^2&0\\ 0&0&(\lambda-1)^3(\lambda+1) \end{matrix} \right].
$$
In Example \ref{examp}, we found that the scalar moment method, i.e. when $\xi=1$, missed the additional multiplicities associated with the polynomial $d_2(\lambda)$. \\
Consider now $\xi=2$ as the size of the block moments $M_k$, the contour $\varphi(t) = 1+\frac{1}{10} e^{\imath t}$, containing 5 eigenvalues $\lambda =1$, as before, and the matrices:
$$
U = \left[ \begin{matrix} 1&0\\ 5&-3\\ 2&-4 \end{matrix} \right], \quad V = \left[ \begin{matrix} 1&3\\ 0&1\\ -2&4 \end{matrix} \right].
$$
We find the block moments $M_k$:
\begin{align*}
&M_0 = \left[ \begin{matrix} -9&-12\\ 9&12 \end{matrix} \right], \quad &M_1 = \left[ \begin{matrix} -1&-22\\ -1&27 \end{matrix} \right], \quad &M_2 = \left[ \begin{matrix} -5&-8\\ 1&18 \end{matrix} \right], \\
&M_3 = \left[ \begin{matrix} -21&30\\ 15&-15 \end{matrix} \right], \quad &M_4 = \left[ \begin{matrix} -49&92\\ 41&-72 \end{matrix} \right], \quad &M_5 = \left[ \begin{matrix} -89&178\\ 79&-153  \end{matrix} \right].
\end{align*}
Then, we have the Hankel matrix $H_{L0}$:
$$
 H_{\xi0} = \left[ \begin{matrix} M_0&M_1&M_2\\ M_1&M_2&M_3\\ M_2&M_3&M_4 \end{matrix} \right] = \left[ \begin{matrix} -9&-12&-1&-22&-5&-8\\ 9&12&-1&27&1&18\\ -1&-22&-5&-8&-21&30\\ -1&27&1&18&15&-15\\ -5&-8&-21&30&-49&92\\ 1&18&15&-15&41&-72 \end{matrix} \right].
$$
The matrix $H_{\xi0}$ is singular. This happens because there are just 5 eigenvalues inside the contour and $H_{\xi0}$ has size $6\times 6$. Then, we have to reduce the matrices $H_{\xi0}$ and $H_{\xi1}$ to match the number of eigenvalues in the contour. Therefore, we get the truncated matrices:
\begin{footnotesize}
$$
 \hat{H}_{\xi0} = \left[ \begin{matrix} -9&-12&-1&-22&-5\\ 9&12&-1&27&1\\ -1&-22&-5&-8&-21\\ -1&27&1&18&15\\ -5&-8&-21&30&-49\end{matrix} \right], \;  \hat{H}_{\xi1} = \left[ \begin{matrix} -1&-22&-5&-8&-21\\ -1&27&1&18&15\\ -5&-8&-21&30&-49\\ 1&18&15&-15&41\\ -21&30&-49&92&-89\end{matrix} \right].
$$
\end{footnotesize}
Then, we obtain:
\begin{equation*}
T = \hat{H}_{\xi0}^{-1}\hat{H}_{\xi1} =  \left[ \begin{matrix} 0&0&0&-2&1\\ 0&0&0&-1&0\\ 1&0&0&4&-3\\ 0&1&0&2&0\\ 0&0&1&-2&3 \end{matrix} \right].
\end{equation*}
The eigenvalues of the matrix $T$ are $1,1,1,1,1$, which are all the eigenvalues inside the contour.\\
Moreover, computing the matrix $Y = [S_0, S_1, S_2]$,  using we get:
$$
\hat{Y} = \left[ \begin{matrix} 0&1&1&2&0\\ 0&-2&-2&0&0\\ 0&-\frac{3}{2}&-\frac{7}{2}&-3&-4 \end{matrix} \right].
$$
Then, $(\hat{Y},T)$ is an invariant pair for $P(\lambda)$.
\end{example}
Experimentally, we noted that the block method allows us to better ``capture'' the multiplicity structure of eigenvalues, when there are several Jordan blocks per eigenvalue. Further investigation of this approach will be the topic of future work. It should be pointed out that the results in \cite{beyn2012integral}, and particularly Theorem 3.3, provide useful insight into a generalized block moment method and into the (good) behavior of the method in presence of multiple eigenvalues. 

Another delicate issue pertaining to contour integral method is 
the choice of $\Gamma$.  If some information about the localization of the eigenvalues is available, one can choose the contour accordingly. In other cases, $\Gamma$ may be taken as a circle for ease of computation, as we do here.

A related question is:  how many eigenvalues of $P(\lambda)$ live inside a given contour? Even an approximate estimate can be useful to choose $\Gamma$ and $k$ consistently. An answer to this problem is provided in \cite{futamura2011stochastic} and \cite{futamura2010parallel}. In particular, Theorem 2 in  \cite{futamura2011stochastic} points out that the number $m$ of eigenvalues of $P(\lambda)$ that are inside $\Gamma$ is given by: 
\begin{equation} \label{ApprEig}
m =  \oint_\Gamma \rm{trace}\left( P(\lambda)^{-1} \frac{dP(\lambda)}{d\lambda} \right) d\lambda.
\end{equation}
For practical computation, the right hand side of equation \eqref{ApprEig} can be approximated by a quadrature rule, thus yielding an estimate for $m$.

Moreover, the choice of $\Gamma$ can be combined with shifting techniques for the eigenvalues of $P(\lambda)$: see for instance \cite{meini2013shift}.

\subsection{Numerical approximation and refinement of invariant pairs}\label{sec_Newton}
When implementing numerical computation of invariant pairs via the (block) moment method, we use numerical quadrature to approximate the moments $\mu_k$ and the vectors $s_k = u^H\mu_k$, respectively defined in \eqref{moments} and \eqref{eq_sk}.  

 \subsubsection{Numerical approximation: trapezoid rule for moments} \label{sec_trapezoid}
Consider the equation \eqref{moments} and assume that $\Gamma$ has a $2 \pi$-periodic smooth parametrization:
\begin{equation*}
\varphi \in C^1(\mathbb{R}, \mathbb{C}), \qquad \varphi(t+2\pi) = \varphi(t) \quad \forall t \in \mathbb{R}.
\end{equation*}
Then, for $k=0,\ldots,2m-1$, we have:
$$
\mu_k = \frac{1}{2\pi  \imath} \oint_\Gamma z^k f(z) dz = \frac{1}{2\pi  \imath} \int_0^{2\pi} \varphi(t)^k f(\varphi(t))\varphi '(t) dt.
$$
Taking equidistant nodes $t_j = \frac{2 j \pi}{N}$, $j=0,\ldots, N-1$, and using the trapezoid rule, we obtain the approximation:
$$
\mu_k \approx \frac{1}{\imath N} \displaystyle \sum_{j=0}^{N-1} \varphi(t_j)^k f(\varphi(t_j))\varphi '(t_j).
$$ 
\subsubsection{Numerical refinement: incorporating line search into Newton's method}\label{sec_LineSearch}
Once an invariant pair has been numerically approximated, it can be refined using an iterative method such as Newton: this is, for instance, the strategy proposed in \cite{betcke2011perturbation}.

Newton's method defines the correction $(\Delta X, \Delta S)$ at each iteration as 
\begin{equation} \label{NM}
P(X,S)+\mathbb{D}P_{(X,S)}(\Delta X, \Delta S)  = 0
\end{equation}

In this section, we show how to incorporate exact line searches into Newton's method for solving the invariant pair problem $P(X,S)=0$. Line searches are relatively inexpensive and improve the global convergence properties of Newton's method (see \cite{wright1999numerical}).\\

\begin{small}

\begin{algorithm} (Newton's Method with Line Search) \label{algor}\\
Input: initial approximation $(X_0,S_0)$, tolerance $\epsilon$.\\ 
Output: better approximation $(X_k,S_k)$ to \eqref{InvPair}. \\
step 1: Set $k=0$\\
step 2: If $\frac{\| P(X_k, S_k) \|_F}{\| X_k \|_F} < \epsilon$: \textbf{STOP} \\
step 3: Solve for $(\Delta X_k,\Delta S_k)$ the equation: 
\begin{equation} \label{correctionEqu}
\mathbb{D}P_{(X,S)}(\Delta X_k, \Delta S_k) = -P(X_k,S_k)
\end{equation}
step 4: Find by exact line searches a $t$ that minimizes the function:
\begin{equation} \label{lineSearchEqu}
\underset{t \in [0,2]}{\min}\| P(X+t\Delta X, S +t\Delta S) \|_F^2
\end{equation}
step 5: Update
\begin{itemize}
\item $X_{k+1} = X_k + t\Delta X_k$, $S_{k+1} = S_k + t\Delta S_k$.
\item $k = k+1$ and go to step 2.
\end{itemize} 
\end{algorithm}
\end{small}
Each iteration of a line search method computes a search direction $d_k$ and then decides how far to move along that direction. The iteration is given by
$$
x_{k+1} = x_{k}+t_k d_k
$$
where the positive scalar $t_k$ is the step length. The success of a line search method depends on effective choices of both the direction $d_k$ and the step length $t_k$ (see \cite{wright1999numerical}). A value $t_k = 1$ gives the original Newton iteration.

In our specific problem \eqref{InvPair}, the direction $d_k$ is given by the solution $(\Delta X_k, \Delta S_k)$ of the correction equation \eqref{correctionEqu}. The step length $t_k$ on each iteration is given by the solution of the minimization problem:
$$
p(t) =  \| P(X+t\Delta X, S +t\Delta S) \|_F^2.
$$
To facilitate the calculation, we consider the following equivalent representation for \eqref{InvPair}. A pair $(X,S) \in \mathbb{C}^{n\times k}\times \mathbb{C}^{k\times k}$ is an invariant pair if and only if satisfies the relation (see \cite{beyn2011continuation}):
\begin{equation} \label{ContInt}
\frac{1}{2 \pi \imath } \oint_\Gamma P(\lambda) X (\lambda I - S)^{-1} d\lambda = 0,
\end{equation}
where $\Gamma \subseteq \mathbb{C}$ is a closed contour with the spectrum of $S$ in its interior. Using the formula for the total derivative of $P$ at $(X,S)$ in direction $(\Delta X, \Delta S)$:
\begin{footnotesize}
\begin{equation*} 
\mathbb{D} P_{(X,S)}(\Delta X, \Delta S) = \frac{1}{2 \pi \imath } \oint_\Gamma P(\lambda) \left(\Delta X + X (\lambda I - S)^{-1} \Delta S\right) (\lambda I - S)^{-1} d\lambda,
\end{equation*}
\end{footnotesize}
we have, at second order in $\|\Delta X\|$ and $\|\Delta S\|$:
\begin{scriptsize}
\begin{align*}
P(X+t\Delta X, S +t\Delta S)  &=  \frac{1}{2 \pi \imath } \oint_\Gamma P(\lambda) (X+t\Delta X) (\lambda I - S - t\Delta S)^{-1} d\lambda = P(X,S) + t \mathbb{D}P_{(X,S)}(\Delta X, \Delta S) +\\ &+ t^2 \left[ \frac{1}{2 \pi \imath } \oint_\Gamma P(\lambda) \left[\Delta X + X (\lambda I - S)^{-1} \Delta S\right] (\lambda I - S)^{-1}\Delta S (\lambda I - S)^{-1} d\lambda \right] +\\
& + t^3 \left[  \frac{1}{2 \pi \imath } \oint_\Gamma P(\lambda) \Delta X (\lambda I - S)^{-1} \Delta S (\lambda I - S)^{-1}\Delta S (\lambda I - S)^{-1} d\lambda  \right] 
\end{align*}
\end{scriptsize}
Recalling that Newton's method defines $(\Delta X, \Delta S)$ by \eqref{NM}, we have:
\begin{scriptsize}
\begin{align*}
P(X+t\Delta X, S +t\Delta S) &=  (1-t)P(X, S) + \\
                                                             & + t^2 \left[ \frac{1}{2 \pi \imath } \oint_\Gamma P(\lambda) \left[\Delta X + X (\lambda I - S)^{-1} \Delta S\right] (\lambda I - S)^{-1}\Delta S (\lambda I - S)^{-1} d\lambda \right] +\\
                                                             & + t^3 \left[  \frac{1}{2 \pi \imath } \oint_\Gamma P(\lambda) \Delta X (\lambda I - S)^{-1} \Delta S (\lambda I - S)^{-1}\Delta S (\lambda I - S)^{-1} d\lambda  \right] 
\end{align*}
\end{scriptsize}
Thus, we have:
\begin{footnotesize}
\begin{align} 
p(t) =& (1-t)^2\| P(X,S)\|_F^2 + t^4\|A\|_F^2 + t^6\|B\|_F^2+ t^2(1-t)\text{trace}( P(X,S)^\ast A + A^\ast P(X,S) ) + \nonumber \\ 
           & +t^3(1-t)\text{trace}(P(X,S)^\ast B + B^\ast P(X,S) ) +t^5\text{trace}(A^\ast B + B^\ast A ) \nonumber\\
\equiv & (1-t)^2\alpha + t^4\theta + t^6\varphi + t^2(1-t)\beta+t^3(1-t)\gamma + t^5\eta \nonumber \\
      =  & t^6\varphi + t^5\eta + t^4(\theta-\gamma) + t^3(\gamma-\beta) + t^2(\alpha+\beta) - 2\alpha t+\alpha
\end{align}
\end{footnotesize}
where:
\begin{footnotesize}
\begin{align*}
& A = \frac{1}{2 \pi \imath } \oint_\Gamma P(\lambda) \left[\Delta X + X (\lambda I - S)^{-1} \Delta S\right] (\lambda I - S)^{-1}\Delta S (\lambda I - S)^{-1} d\lambda,\\
& B = \frac{1}{2 \pi \imath } \oint_\Gamma P(\lambda) \Delta X (\lambda I - S)^{-1} \Delta S (\lambda I - S)^{-1}\Delta S (\lambda I - S)^{-1} d\lambda,\\
& \alpha = \| P(X,S)\|_F^2, \; \theta = \|A\|_F^2, \; \varphi  = \|B\|_F^2, \eta = \text{trace}(A^\ast B + B^\ast A),\\
& \beta =  \text{trace}(P(X,S)^\ast A + A^\ast P(X,S) ), \gamma =\text{trace}(P(X,S)^\ast  B + B^\ast P(X,S) ).
\end{align*}
\end{footnotesize}
Therefore, in each iteration of Algorithm \ref{algor}, solving \eqref{lineSearchEqu} is equivalent to finding the minimum of the polynomial $p(t)$ for $t \in [0,2]$. 

\subsection{Numerical results}
In this section we compare two methods to refine approximate invariant pairs $(X,S) \in \mathbb{C}^{n \times k}\times \mathbb{C}^{k\times k}$: Newton's method (N.M.) presented in \cite{betcke2011perturbation} and Newton's method with line search (N.M.L.S.), explained in Section \ref{sec_LineSearch}.

We have implemented both methods in MATLAB and applied them to several problems taken from the NLEVP collection (see \cite{betcke2013nlevp}). For each problem, an initial invariant pair $(X_0, S_0)$ has first been approximated using the (block) moment method of Section \ref{sec_computingInvPairs} and approximating the moments $\mu_i$ in \eqref{moments} via the trapezoid rule discussed in Section \ref{sec_trapezoid}, with $N=20$ integration nodes. Moreover, $\Gamma$ is chosen for each problem as the contour enclosing the $k$ eigenvalues with largest condition number (computed using the MATLAB function \texttt{polyeig}).

Table \ref{table_Newton} shows that line search is generally effective in reducing the number of iterations and the overall computation time.
\begin{table}[h!]
\centering
\begin{footnotesize}
\begin{tabular}{|c|c|c||c|c||c|c|c|}
 \multicolumn{5}{r}{N.M.}&\multicolumn{2}{r}{N.M.L.S.}\\
 \hline \textbf{Problem} & \textbf{Deg $P(\lambda)$} &  \textbf{$n\times k$} & \textbf{Ite} & \textbf{Time} & \textbf{Ite} & \textbf{Time}\\
\hline bicycle&2&$2\times 2$&23&0.082&16&0.112\\
\hline butterfly&4&$64\times 5$&67&3.719&22&1.567\\ 
\hline cd\_player&2&$60\times 6$&500&N.C.&19&1.021\\
\hline closed\_loop&2&$2\times 2$&8&0.016&7&0.02\\
\hline damped\_beam&2&$200\times 6$&28&6.109&4& 0.6037\\
\hline dirac&2&$80\times 6$&500&N.C.&41&1.965\\   
\hline hospital&2&$24\times 24$&53&5.65&51&6.21\\ 
\hline metal\_strip&2&$9\times 9$&500&N.C.&28&0.589\\   
\hline mobile\_manipulator&2&$5\times 2$&8&0.014&7&0.030\\
\hline pdde\_stability&2&$225\times 6$&29&9.644&16&5.622\\ 
\hline planar\_waveguide&4&$129\times 6$&72&11.148&19&3.682\\
\hline plasma\_drift&3&$128\times 6$&69&13.059&26&5.596\\
\hline power\_plant&2&$8\times 8$&15&0.34&13&0.39\\
\hline railtrack&2&$1005\times 3$&32&199.365&28&209.471\\
\hline
\end{tabular}
\end{footnotesize}
\caption{Comparison of results for classical Newton and Newton with line search.}\label{table_Newton} \vspace{-0.5em}
\end{table}

\begin{figure}[h!]
 \begin{center}
  \includegraphics[height=17em]{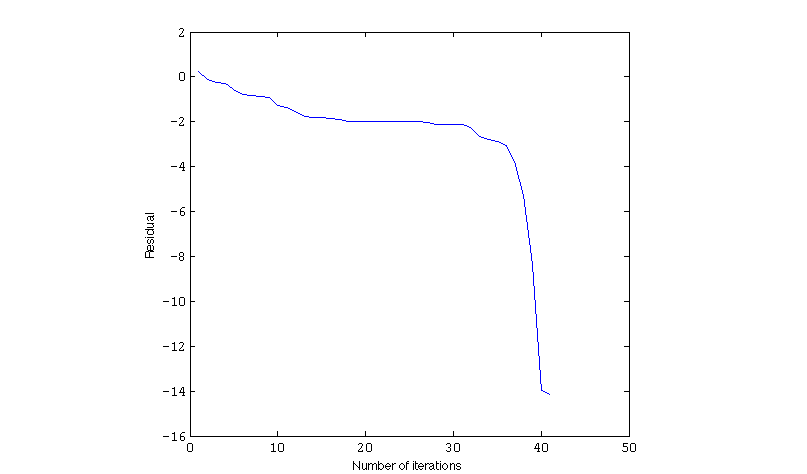}
\caption{Convergence of Dirac problem using Newton's method with line search.This is a log-10 plot of the relative residual $\frac{\|P(X,S)\|_F}{\|X\|_F}$ versus the number of iterations.}\label{figure} \vspace{-2em}
  \end{center}
  \end{figure}

%
Figure \ref{figure} shows the convergence of the Newton's method with line search, for the Dirac problem presented in Table \ref{table_Newton}. Here we use as contour the circle of center $C = -0.1$ and radius $R = 1.14$, which contains the 6 eigenvalues with largest condition number.

\section{Matrix solvents}\label{sec_solvents}
In this section we study the matrix solvent problem as a particular case of the invariant pair problem, and we apply to solvents some results that we have obtained for invariant pairs.
\begin{definition}
A matrix $S \in \mathbb{C}^{n\times n}$ is called a solvent for $P(S)$ if satisfies the relation:
\begin{equation} \label{Sol}
P(S):=A_\ell S^\ell +\cdots+A_2S^2 + A_1S + A_0 = 0.
\end{equation}
\end{definition}
A special case is, for $\ell =2$, the quadratic matrix equation $Q(S) := A_2S^2 +A_1S+A_0 = 0$, which has received considerable attention in the literature. For instance, in \cite{higham2000numerical} and \cite{higham2001solving} the authors find formulations for the condition number and the backward error. They also propose functional iteration approaches based on Bernoulli's method and Newton's method with line search to compute the solution numerically.

The relation between eigenvalues of $P(\lambda)$ and solvents is highlighted in \cite{lancaster2002lambda}: a corollary of the generalized B\'ezout theorem states that if $S$ is a solvent of $P(S)$ then:
$$
P(\lambda) = L(\lambda)(\lambda I - S),
$$
where $L(\lambda)$ is a matrix polynomial of degree $\ell -1$. Then any eigenpair of the solvent $S$ is an eigenpair of $P(\lambda)$.

\subsection{Condition number and backward error}
An analysis and a computable formulation for the condition number and backward error of the quadratic matrix equation $Q(S)=0$ can be found in \citep{higham2000numerical} and \citep{higham2001solving}.

Here we give explicit expressions for the condition number and backward error for the general matrix solvent problem $P(S)=0$. We  follow the ideas presented in \citep{higham2000numerical}, \citep{higham2001solving} and \citep{tisseur2000backward}.
\subsubsection{Condition number} \label{secCondNumSol}
We perform here a similar analysis as we did in Section \ref{SecCondNumInvPa}. \\
A normwise condition number of the solvent $S$ can be defined by:
\begin{footnotesize}
\begin{equation} \label{defCondNumSolv}
\kappa(S) = \limsup_{\epsilon\rightarrow 0} \left\lbrace \frac{1}{\epsilon} \frac{\left|\left| \Delta S \right|\right|_F}{\left|\left| S \right|\right|_F} : ( P+\Delta P)(S+\Delta S)=0, \right. \| \Delta A_i\|_F \leq \epsilon \alpha_i, i=0:\ell \Bigg\},
\end{equation}
\end{footnotesize}
where $\Delta P(\lambda) = \displaystyle \sum_{i=0}^\ell \lambda^i \Delta A_i$. The $\alpha_i$ are nonnegative weights; in particular, $\Delta A_i$ can be forced to zero by setting $\alpha_i =0$.
\begin{theorem} \label{condNumSolvents}
The normwise condition number of the solvent $S$ is given by:
$$
\kappa(S) = \frac{\left|\left| \hat{B}_{S}^{-1} \hat{B}_{A}\right|\right|_2}{\left|\left| S \right|\right|_F},
$$
where
\begin{footnotesize}
$$
\hat{B}_{S} = \displaystyle\sum_{j = 1}^{\ell}  \displaystyle\sum_{i=0}^{j-1}\left( (S^{j-i-1})^T\otimes A_j S^i \right), \: \hat{B}_{A} = \left[ \begin{matrix} \alpha_\ell (S^\ell)^T\otimes I_n & \alpha_{\ell-1}(S^{\ell-1})^T\otimes I_n & \cdots &  \alpha_0 I_{n^2}  \end{matrix} \right]
$$
\end{footnotesize}
\end{theorem}
The proof of Theorem \ref{condNumSolvents} follows from the proof of Theorem \ref{condNumInvPair}, by taking $\Delta X = 0$, $X=I$ and noting that the matrix $S$ has size $n\times n$ in the matrix solvent problem.

\subsubsection{Backward error}
Let $\alpha_i$, for $i=0,\ldots, \ell$, be nonnegative weights as in Section \ref{secCondNumSol}. The backward error of an approximate solution $T$ to \eqref{Sol} can be defined as:
\begin{equation} \label{BES1}
\eta(T)  = \min \{\epsilon: (P+\Delta P)(T) = 0, \| \Delta A_i\|_F \leq \epsilon \alpha_i, i=0,\ldots, \ell \}
\end{equation}
We proceed as in Section \ref{backErrorInv} and we obtain bounds for the backward error of $P(S)$: 
\begin{equation*}
\eta(T) \geq \frac{\| P(T) \|_F}{(\alpha_\ell^2 \| T^\ell \|_F^2 + \cdots + \alpha_1^2 \| T \|_F^2 + \alpha_0^2)^{1/2}},
\end{equation*}
\begin{equation*}
\eta(T) \leq \frac{\| P(T) \|_F}{(\alpha_\ell^2\sigma_{\min}(T^\ell)^2 + \cdots + \alpha_1^2\sigma_{\min}(T)^2 + \alpha_0^2)^{1/2}},
\end{equation*}
where $\sigma_{\min}$ denotes the smallest singular value, that by assumption is nonzero.

\section{Computation of solvents}\label{sec_solvents_computation}

Motivated by applications to differential equations \cite{brugnano1998solving}, we study an approach to the computation of solvents based on the moment method, by specializing the results presented in Section \ref{sec_invpairs}.

Let us recall some results that will be needed later. The next result is a generalization of a theorem presented in \cite{higham2001solving} which gives information about the number of solvents of $P(S)$.
\begin{theorem} \label{Teo1}
Suppose $P(\lambda)$ has $p$ distinct eigenvalues $\left\lbrace\lambda_i\right\rbrace^p_{i=1}$, with $n\leq p \leq \ell n$, and that the corresponding set of $p$ eigenvectors $\left\lbrace v_i \right\rbrace^p_{i=1}$ satisfies the Haar condition (every subset of $n$ of them is linearly independent). Then there are at least $\left( \begin{matrix} p\\n\end{matrix} \right)$ different solvents of $P(\lambda)$, and exactly this many if $p=\ell n$, which are given by 
$$
S = W \text{diag}(\mu_i)W^{-1}, \quad W = \left[ \begin{matrix} w_1 & \cdots & w_n\end{matrix} \right],
$$
where the eigenpairs $(\mu_i,w_i)_{i=1}^n$ are chosen from among the eigenpairs $(\lambda_i,v_i)_{i=1}^p$ of $P$.
\end{theorem}
Note that if we have that $p=n$ in Theorem \ref{Teo1}, the distinctness of the eigenvalues is not needed, and then we have a sufficient condition for the existence of a solvent.
\begin{corollary}
If $P(\lambda)$ has $n$ linearly independent eigenvectors $v_1,v_2,\ldots,v_n$ then $P(S)$ has a solvent. 
\end{corollary}  

An example which illustrates this last result is the following. Consider the quadratic matrix solvent problem (see \citep{dennis1976algebraic}, \cite{higham2001solving})
$$
Q(S) = S^2 +  \left[ \begin{matrix} -1 &-6\\2 &-9 \end{matrix} \right] S + \left[ \begin{matrix} 0 &12\\-2 &14 \end{matrix} \right]
$$
$Q(\lambda)$ has eigenpairs: $\left( 1,\left[ \begin{matrix} 1\\0 \end{matrix} \right] \right)$, $\left( 2,\left[ \begin{matrix} 0\\1 \end{matrix} \right] \right)$, $\left( 3,\left[ \begin{matrix} 1\\1 \end{matrix} \right] \right)$ and $\left( 4,\left[ \begin{matrix} 1\\1 \end{matrix} \right] \right).$
For this example, the complete set of solvents is:
$$
\left[ \begin{matrix} 1&0\\0&2 \end{matrix} \right], \left[ \begin{matrix} 1&2\\0&3 \end{matrix} \right], \left[ \begin{matrix} 3&0\\1&2 \end{matrix} \right], \left[ \begin{matrix} 1&3\\0&4 \end{matrix} \right] \text{and} \left[ \begin{matrix} 4&0\\2&2 \end{matrix} \right].
$$
Note that we cannot construct a solvent whose eigenvalues are 3 and 4 because the associated eigenvectors are linearly dependent.

Our approach to compute matrix solvents is based on the relation between the matrix solvent problem \eqref{Sol} and the invariant pair problem \eqref{InvPair}. We state this in the following result.
\begin{theorem}
Let $P(\lambda)$ be a $n\times n$ matrix polynomial and consider an invariant pair $(Y,T) \in \mathbb{C}^{n\times k}\times \mathbb{C}^{k\times k}$ of $P(\lambda)$. If the matrix $Y$ has size $n\times n$, i.e. $k=n$, and is invertible, then $S=YTY^{-1}$ satisfies equation \eqref{Sol}, i.e., $S$ is a matrix solvent of $P(\lambda)$.
\begin{proof}
As $(Y,T) \in \mathbb{C}^{n\times n}\times \mathbb{C}^{n\times n}$ is an invariant pair of $P(\lambda)$, we have:
$$
A_\ell YT^\ell +\cdots+A_2YT^2 + A_1YT + A_0Y = 0.
$$
Since $Y$ is invertible, we can post-multiply by $Y^{-1}$. Then we get:
\begin{align*}
& A_\ell YT^\ell Y^{-1} +\cdots+A_2YT^2Y^{-1} + A_1YTY^{-1} + A_0 = 0 \Leftrightarrow \\
& A_\ell S^\ell +\cdots+A_2S^2 + A_1S + A_0 = 0.
\end{align*}
Therefore, $S$ is a matrix solvent of $P(\lambda)$. 
\end{proof}
\end{theorem}


\subsection{Numerical refinement of solvents}
As pointed out before, the use of Newton's method incorporating line searches to find solvents is not new. For instance, in \cite{higham2001solving}, \cite{long2008improved} this approach is used to approximate solvents for the quadratic matrix equation. Here we apply this method to approximate solvents for the general matrix solvent problem $P(S)=0$ and follow the ideas of Section \ref{sec_LineSearch}.

The application of Newton's method with line search to find solvents is based on the following steps for the $k$-th iteration.
\begin{itemize} 
\item Solve for $\Delta S_k$ the equation: $\mathbb{D}P_{S}(\Delta S_k) = -P(S_k)$
\item Find by exact line searches a $t_k$ that minimizes the function:
$$
\underset{t_k \in [0,2]}{\min}\| P(S_k +t_k\Delta S_k) \|_F^2 
$$
\item Update $S_{k+1} = S_k + t_k \Delta S_k$.
\end{itemize}
The step length $t_k$ on each iteration is given by the solution of the minimization problem:
$$
p(t_k) =  \| P(S_k +t_k\Delta S_k) \|_F^2,
$$
As we did in Section \ref{sec_LineSearch} for the invariant pair problem, we use an equivalent contour integral representation for \eqref{Sol}. A matrix $S \in \mathbb{C}^{n\times n}$ is a solvent if and only if satisfies the relation:
$$
\frac{1}{2 \pi \imath } \oint_\Gamma P(\lambda) (\lambda I - S)^{-1} d\lambda = 0,
$$
for any closed contour $\Gamma \subseteq \mathbb{C}$ with the spectrum of $S$ in its interior.\\
Using the formula for the total derivative of $P$ at $S$ in direction $\Delta S$:
\begin{equation} \label{derCont}
\mathbb{D}P_{S}(\Delta S) = \frac{1}{2 \pi \imath } \oint_\Gamma P(\lambda) (\lambda I - S)^{-1} \Delta S (\lambda I - S)^{-1} d\lambda,
\end{equation}
we obtain:
\begin{scriptsize}
$$
P(S +t\Delta S)  = P(S) + t \mathbb{D}P_S(\Delta S) + t^2 \left[ \frac{1}{2 \pi \imath } \oint_\Gamma P(\lambda) (\lambda I - S)^{-1} \Delta S (\lambda I - S)^{-1}\Delta S (\lambda I - S)^{-1} d\lambda \right].
$$
\end{scriptsize}
Recalling that Newton's method defines $\Delta S$ by 
\begin{equation*}
P(S)+\mathbb{D}P_S (\Delta S)  = 0
\end{equation*}
then we have:
\begin{footnotesize}
$$
P(S +t\Delta S) =  (1-t)P(S) + t^2 \left[ \frac{1}{2 \pi \imath } \oint_\Gamma P(\lambda) (\lambda I - S)^{-1} \Delta S (\lambda I - S)^{-1}\Delta S (\lambda I - S)^{-1} d\lambda \right].
$$
\end{footnotesize}
Thus, we obtain:
\begin{align*} 
p(t) =  &\: (1-t)^2\| P(S)\|_F^2 + t^4\|A\|_F^2 + t^2(1-t)\text{trace}(P(S)^\ast A + A^\ast P(S) ) \nonumber\\
\equiv &\: (1-t)^2\alpha + t^4\theta + t^2(1-t)\beta = \nonumber\\
      =  &\: t^4\theta - t^3\beta + t^2(\alpha+\beta) - 2\alpha t+\alpha
\end{align*}
where:
\begin{align*}
A & = \frac{1}{2 \pi \imath } \oint_\Gamma P(\lambda) (\lambda I - S)^{-1} \Delta S (\lambda I - S)^{-1}\Delta S (\lambda I - S)^{-1} d\lambda,\\
\alpha &= \|P(S) \|_F^2, \quad \theta = \|A\|_F^2, \quad \beta =  \text{trace}(P(S)^\ast A + A^\ast P(S) ).
\end{align*}
Therefore, in each iteration one finds the minimum of the polynomial $p(t)$ for $t \in [0,2]$. 

\section{Solvents and triangularized matrix polynomials}\label{sec_solvents_triang}
Motivated by the results in \cite{taslaman2013triangularizing} and \cite{tisseur2013triangularizing}, where the authors analyze a method for triangularizing the matrix polynomial $P(\lambda)$, we aim here to study the relation between solvents of general and of triangularized matrix polynomials.

\subsection{Triangularizing matrix polynomials}
For any algebraically closed field $\mathbb{F}$, any matrix polynomial $P(\lambda) \in \mathbb{F}[\lambda]^{n\times m}$, with $n\leq m$, can be reduced to triangular form via unimodular transformations, preserving the degree and the finite and infinite elementary divisors \cite{taslaman2013triangularizing}, \cite{tisseur2013triangularizing}.

\begin{theorem} \cite{taslaman2013triangularizing}
For an algebraically closed field $\mathbb{F}$, any $P(\lambda) \in \mathbb{F}[\lambda]^{n\times m}$ with $n\leq m$ is triangularizable. 
\end{theorem}

What can we say about solvents for a given matrix polynomial $P(\lambda)$ and for the associated triangularized polynomial? A partial answer will be given in Theorem \ref{thm_triang}.

 \begin{theorem} \label{T3}
For any $\ell n\times \ell n$ monic linearization $\lambda I - A$ of $P(\lambda) \in \mathbb{C}[\lambda]^{n\times n}$ with nonsingular leading coefficient, there exists $U \in \mathbb{C}^{n\times \ell n}$ with orthogonal columns such that $M = \left[ \begin{matrix} U\\ UA\\ \vdots\\ UA^{\ell-1} \end{matrix} \right]$ is nonsingular and $\lambda I - M A M^{-1} $ is a linearization for the polynomial
$T(\lambda) = \lambda^\ell I + \lambda^{\ell-1} T_{\ell-1} + \cdots+\lambda^2 T_2 + \lambda T_1 + T_0$, which is upper triangular and equivalent to $P(\lambda)$.
\end{theorem}
Theorem \ref{T3} is a straightforward generalization of a result found in \cite{tisseur2013triangularizing}. Note that, for the time being, we assume that the leading coefficient $A_\ell$ is nonsingular.
\begin{theorem} \label{thm_triang}
Let $P(\lambda)$ be a $n\times n$ matrix polynomial and consider the linearization:
$$
A =  \left[ \begin{matrix} 0&I_n&0&\cdots&0\\ 0&0&I_n&\cdots&0\\ \vdots&\vdots&\vdots&\ddots&\vdots\\ 0&0&0&\cdots&I_n\\ -A_0&-A_1&-A_2&\cdots&-A_{\ell-1}\\ \end{matrix} \right].
$$
Let $M$ be as in Theorem \ref{T3} and let $Y_1$ be the first $n\times n$ block of the matrix \\$M^{-1} \left[ \begin{matrix} I_n\\ S_t\\ \vdots \\ S_t^{\ell-1} \end{matrix} \right]$. If $Y_1$ is nonsingular and $S_t$ is a solvent for the triangularized problem, i.e., $T(S_t) = 0$, then $S = Y_1S_tY_1^{-1}$ is a solvent for $P(S)$.
\begin{proof}
Note that:
\begin{footnotesize}
\begin{align*}
\left[ \begin{matrix}  0\\ 0\\ \vdots\\ 0 \end{matrix} \right] &= \left[ \begin{matrix}  S_t - S_t\\S_t^2 - S_t^2\\ \vdots\\  -T_0-T_1S_t - T_2S^2_t-\cdots-S_t^{\ell} \end{matrix} \right] =\\
& = \left[ \begin{matrix} 0&I_n&0&\cdots&0\\ 0&0&I_n&\cdots&0\\ \vdots&\vdots&\vdots&\ddots&\vdots\\ 0&0&0&\cdots&I_n\\ -T_0&-T_1&-T_2&\cdots&-T_{\ell-1}\\ \end{matrix} \right]  \left[ \begin{matrix}  I_n\\ S_t\\ \vdots\\ S_t^{\ell-1} \end{matrix} \right]- \left[ \begin{matrix}  S_t \\ S^2_t\\ \vdots\\ S^{\ell}_t\end{matrix} \right] = \\
&= M^{-1} \left[ \begin{matrix} 0&I_n&0&\cdots&0\\ 0&0&I_n&\cdots&0\\ \vdots&\vdots&\vdots&\ddots&\vdots\\ 0&0&0&\cdots&I_n\\ -T_0&-T_1&-T_2&\cdots&-T_{\ell-1}\\ \end{matrix} \right] M M^{-1} \left[ \begin{matrix}  I_n\\ S_t\\ \vdots\\ S_t^{\ell-1} \end{matrix} \right] - M^{-1} \left[ \begin{matrix}  I_n\\ S_t\\ \vdots\\ S_t^{\ell-1} \end{matrix} \right] S_t \stackrel{(iii)}{=} \\
& = \left[ \begin{matrix} 0&I_n&0&\cdots&0\\ 0&0&I_n&\cdots&0\\ \vdots&\vdots&\vdots&\ddots&\vdots\\ 0&0&0&\cdots&I_n\\ -A_0&-A_1&-A_2&\cdots&-A_{\ell-1}\\ \end{matrix} \right] M^{-1} \left[ \begin{matrix}  I_n\\ S_t\\ \vdots\\ S_t^{\ell-1} \end{matrix} \right] - M^{-1} \left[ \begin{matrix}  I_n\\ S_t\\ \vdots\\ S_t^{\ell-1} \end{matrix} \right] S_t.
\end{align*}
\end{footnotesize}
Since $M^{-1} \left[ \begin{matrix}  I_n\\ S_t\\ \vdots\\ S_t^{\ell-1} \end{matrix} \right]$ has size $\ell n\times n$, let us partition it as $\left[ \begin{matrix}  Y_1\\ Y_2\\ \vdots\\ Y_\ell \end{matrix} \right]$, where $Y_i \in \mathbb{C}^{n\times n}$ for $i=1,\ldots,\ell$. Then:
\begin{equation*}
\left[ \begin{matrix} 0&I_n&0&\cdots&0\\ 0&0&I_n&\cdots&0\\ \vdots&\vdots&\vdots&\ddots&\vdots\\ 0&0&0&\cdots&I_n\\ -A_0&-A_1&-A_2&\cdots&-A_{\ell-1}\\ \end{matrix} \right] \left[ \begin{matrix}  Y_1\\ Y_2\\ \vdots\\ Y_\ell \end{matrix} \right] - \left[ \begin{matrix}  Y_1\\ Y_2\\ \vdots\\ Y_\ell \end{matrix} \right] S_t = \left[ \begin{matrix}  0 \\ 0 \end{matrix} \right].
\end{equation*}
Then we have: 
\begin{align}
&Y_i = Y_1S_t^{i-1}, \; \text{for}\: i=2,\ldots,\ell;  \label{EQ1}\\
\quad -&A_0Y_1 - A_1Y_2-\cdots- A_{\ell-1}Y_\ell - Y_\ell S_t= 0. \label{EQ2}
\end{align}
Substituting equations \eqref{EQ1} in \eqref{EQ2} we obtain:
\begin{equation*}
0 = Y_1S^\ell_t + A_{\ell-1} Y_1S_t^{\ell-1} +\cdots+ A_1Y_1S_t + A_0Y_1.
\end{equation*}
If $Y_1$ is invertible we have:
\begin{equation*}
0 = Y_1S^\ell_t Y_1^{-1} + A_{\ell-1} Y_1S_t^{\ell-1} Y_1^{-1} +\cdots+ A_1Y_1S_t Y_1^{-1} + A_0.
\end{equation*}
Taking $S = Y_1S_tY_1^{-1}$ we have:
\begin{equation*}
0 = S^\ell +A_{\ell-1}S^{\ell-1} + \cdots +A_2S^2 + A_1S  + A_0 := P(S).
\end{equation*} 
Then $S = Y_1S_tY_1^{-1}$ is a solvent for $P(S)$.
\end{proof}
\end{theorem}

\subsection{Example: A problem with an infinite number of solvents}
What happens to the ideas outlined above when working on problems with an infinite number of solvents? Here is an example taken from \cite{pereira2003solvents}. 

Consider the matrix polynomial:
\begin{equation*}
P(\lambda) = \lambda^2 I + \lambda \left[ \begin{matrix} -7&-2&-2\\ \frac{3}{31}&\frac{-203}{31}&\frac{8}{31}\\ \frac{-13}{31}&\frac{-40}{31}&\frac{-231}{31} \end{matrix} \right] + \left[ \begin{matrix} 13&9&7\\ \frac{-21}{31}&\frac{294}{31}&\frac{-36}{31}\\ \frac{60}{31} & \frac{183}{31}&\frac{435}{31} \end{matrix} \right] 
\end{equation*}
Triangularizing $T(\lambda)$ we find: 
\begin{align*}
T(\lambda) &= \left[ \begin{matrix} (\lambda-3)(\lambda-4) & (\lambda-3) & 0\\ 0 & (\lambda-3)^2 & 1\\ 0 & 0 & (\lambda-4)^2 \end{matrix} \right] = \lambda^2 I_2 + \lambda T_1 + T_0 = \\
&= \lambda^2 I + \lambda \left[ \begin{matrix} -7&1&0\\ 0&-6&0\\ 0&0&-8 \end{matrix} \right] + \left[ \begin{matrix} 12&-3&0\\ 0&9&1\\ 0&0&16 \end{matrix} \right].
\end{align*}
Now, suppose that the solvent $S_t \in \mathbb{C}^{3\times 3}$ of the triangularized problem is in upper triangular form, i.e.:
\begin{equation*}
S_t = \left[ \begin{matrix}  x_{11}&x_{12}&x_{13}\\ 0&x_{22}&x_{23}\\ 0&0&x_{33} \end{matrix} \right],
\end{equation*}
then we have:
\begin{tiny}
\begin{align*}
T(S_t) =&  S_t^2 + T_1 S_t + T_0 = \\
=&\left[ \begin{matrix}  (x_{11}-3)(x_{11}-4) & x_{22}-7x_{12}+x_{11}x_{12}+x_{12}x_{22}-3 & x_{23}-7x_{13}+x_{11}x_{13}+x_{12}x_{23}+x_{13}x_{33}\\ 0 & (x_{22}-3)^2 & x_{22}x_{23}-6x_{23}+x_{23}x_{33}+1\\ 0&0&(x_{33}-4)^2 \end{matrix} \right]
\end{align*}
\end{tiny}
In the task of solving the problem $T(S_t) = 0$, we see that: $x_{11} =3$ or $x_{11}=4$, $x_{22}=3$ and $x_{33}=4$. Then we have two cases:
\begin{enumerate}
\item[I.] If $x_{11} =3$, $x_{22} =3$ and $x_{33} =4$:\\
Then we find that $x_{23} =-1$, $x_{12} =0$ and $x_{12} =-1$, which is a contradiction. In this case there is no solution and then we can't construct a solvent.
\item[II.] If $x_{11} =4$, $x_{22} =3$ and $x_{33} =4$:\\
Then we find that $x_{23} =-1$ and $x_{13} =x_{12}+1$. In this case the solvent $S_t$ has the form:
\begin{equation*}
 S_t = \left[ \begin{matrix} 4&x_{12}&x_{12}+1\\ 0&3&-1\\ 0&0&4 \end{matrix} \right] = \left[ \begin{matrix} 4&0&1\\ 0&3&-1\\ 0&0&4 \end{matrix} \right] + x_{12} \left[ \begin{matrix} 0&1&1\\ 0&0&0\\ 0&0&0 \end{matrix} \right], 
\end{equation*}
for $x_{12}\in \mathbb{C}$.
\end{enumerate}
Thus $T(\lambda)$ has an infinite number of solvents and the same holds for $P(\lambda)$.

\section{Conclusion}

In this paper we have explored several questions related to invariant pairs and solvents of matrix polynomials. In particular, preliminary results on the use of scalar or block moment Hankel pencils to compute invariant pairs and solvents suggest that this approach may present several points of interest. A more detailed analysis, along with the design and development of effective algorithms and extensive numerical tests, will be the topic of future work. 

\vspace{10pt}

{\bf Acknowledgments. }We would like to thank Fran\c coise Tisseur for interesting and fruitful discussions, and Daniel Kressner for providing the MATLAB implementation of the algorithm presented in \cite{betcke2011perturbation}. Our thanks also go to two anonymous referees for their suggestions and help in improving this work.


\bibliographystyle{abbrv}
\bibliography{ourpaper}{}

  
\end{document}